\pdfoutput=1
\documentclass{article}

\usepackage{libertine}
\usepackage[T2A,T1]{fontenc}
\usepackage[utf8]{inputenc}      \usepackage[libertine]{newtxmath}
\usepackage[scaled=0.83]{beramono}
\usepackage{eucal}
\usepackage[babel]{microtype}
\usepackage[russian,main=american]{babel}
\usepackage[strict]{csquotes}
\usepackage{alphabeta}              \usepackage[fleqn,leqno]{mathtools} \usepackage[fleqn,leqno]{amsmath}   \usepackage[acronym,toc,nomain]{glossaries}  \usepackage[allcolors=.]{hyperref}               \usepackage[capitalise]{cleveref}   \usepackage[dvipsnames]{xcolor}     \usepackage[ocgcolorlinks]{ocgx2}[2017/03/30]	

                                             \let\digamma\relax                                                     \usepackage{amsfonts}
\usepackage{amsthm}
\usepackage{booktabs}
\usepackage{changepage}
\usepackage{diagbox}
\usepackage{dynkin-diagrams}
\usepackage{enumitem}
\usepackage{parskip}
\usepackage{pgfplots}
\pgfplotsset{compat=1.16}     \usepackage{rotating}
\usepackage{siunitx}
\usepackage{subcaption}
\usepackage{thmtools}
\usepackage{tikz}
\usepackage{tikz-cd}
\usetikzlibrary{babel}
\usepackage[colorinlistoftodos, textsize = footnotesize]{todonotes}
\usepackage{xparse}
\usepackage{xspace}
\usetikzlibrary{calc,tqft,decorations.markings}

\usepackage[backend=biber,maxbibnames=99,doi=false]{biblatex}

\newtoggle{nome}
\toggletrue{nome}

\mathtoolsset{centercolon}

\newcounter{todocounter}
\DeclareDocumentCommand\addreference{g}{\stepcounter{todocounter}\todo[color = blue!30]{\thetodocounter. Add reference\IfNoValueF{#1}{: #1}}\xspace}
\DeclareDocumentCommand\checkthis{g}{\stepcounter{todocounter}\todo[color = red!50]{\thetodocounter. Check this\IfNoValueF{#1}{: #1}}\xspace}
\DeclareDocumentCommand\fixthis{g}{\stepcounter{todocounter}\todo[color = orange!50]{\thetodocounter. Fix this\IfNoValueF{#1}{: #1}}\xspace}
\DeclareDocumentCommand\expand{g}{\stepcounter{todocounter}\todo[color = green!50]{\thetodocounter. Expand\IfNoValueF{#1}{: #1}}\xspace}

\tikzset{every loop/.style={min distance=1cm, looseness=30}}
\tikzset{vertex/.style={circle, draw, black, minimum size=6pt, inner sep=0pt}}
\tikzset{half-vertex/.style={semicircle, draw, fill, black, minimum size=3pt, inner sep=0pt, yshift=1.5pt}}

\tikzset{middlearrow/.style={decoration={markings, mark= at position 0.5 with {\arrow{#1}},}, postaction={decorate}}}

\declaretheoremstyle[
  spaceabove = 3pt,
  spacebelow = 3pt,
  bodyfont = \itshape,
]{first}
\declaretheorem[numberwithin=section, style=first]{theorem}

\declaretheorem[sibling=theorem, style=first]{corollary}
\declaretheorem[sibling=theorem, style=first]{lemma}
\declaretheorem[sibling=theorem, style=first]{proposition}

\declaretheorem[numberwithin=section, style=first, title=Theorem]{alphatheorem}

\declaretheoremstyle[
  spaceabove = 0pt,
  spacebelow = 0pt,
]{second}
\theoremstyle{second}
\declaretheorem[style=second, sibling=theorem]{construction}
\declaretheorem[style=second, sibling=theorem]{definition}
\declaretheorem[style=second, sibling=theorem]{example}
\declaretheorem[style=second, sibling=theorem]{notation}

\declaretheorem[style=second, sibling=theorem]{remark}

\crefname{conjecture}{Conjecture}{Conjectures}
\crefname{construction}{Construction}{Constructions}

\crefname{alphatheorem}{Theorem}{Theorems}
\crefname{alphaconjecture}{Conjecture}{Conjectures}
\crefname{alphacorollary}{Corollary}{Corollaries}
\crefname{alphaproposition}{Proposition}{Propositions}

\crefname{figure}{Figure}{Figures}
\crefname{question}{Question}{Questions}

\makeatletter
\def\gitfootnote{\gdef\@thefnmark{}\@footnotetext}
\makeatother

\makeatletter
\DeclareRobustCommand\widecheck[1]{{\mathpalette\@widecheck{#1}}}
\def\@widecheck#1#2{\setbox\z@\hbox{\m@th$#1#2$}\setbox\tw@\hbox{\m@th$#1\widehat{\vrule\@width\z@\@height\ht\z@
          \vrule\@height\z@\@width\wd\z@}$}\dp\tw@-\ht\z@
    \@tempdima\ht\z@ \advance\@tempdima2\ht\tw@ \divide\@tempdima\thr@@
    \setbox\tw@\hbox{\raise\@tempdima\hbox{\scalebox{1}[-1]{\lower\@tempdima\box
\tw@}}}{\ooalign{\box\tw@ \cr \box\z@}}}
\makeatother

\makeindex

\crefformat{enumi}{#2\textup{(#1)}#3}

\makeatletter
\newcommand{\acro}[3]{\newacronym{#1}{#2}{#3}\@namedef{#1}{{\acrshort{#1}}\xspace}\@namedef{#2}{{\acrlong{#1}}\xspace}}

\acro{qft}{QFT}{Quantum Field Theory}
\acro{tqft}{TQFT}{Topological Quantum Field Theory}
\acro{cft}{CFT}{Conformal Field Theory}
\acro{rcft}{RCFT}{Rational Conformal Field Theory}
\acro{sft}{SFT}{Symplectic Field Theory}
\acro{cohft}{CohFT}{Cohomological Field Theory}
\acro{trft}{TropFT}{Tropical Field Theory}

\acro{syz}{SYZ}{Strominger--Yau--Zaslow}
\acro{tuy}{TUY}{Tsuchiya--Ueno--Yamada}
\acro{wdvv}{WDVV}{Witten--Dijkgraaf--Verlinde--Verlinde}
\acro{wz}{WZ} {Wess--Zumino--Witten}

\acro{glsm}{GLSM}{Gauged Linear Sigma Model}
 \newcommand   \graph {\gamma}          \newcommand \sphere  {\ensuremath{S}}                    \newcommand\LT{\ensuremath{\mathrm{L}^2(\mathrm{S}^1;\mathbb{C})}}

\newcommand\ii{\ensuremath{\sqrt{-1}}}
\newcommand \dlog [1] {\frac{d#1}{#1}}

                                                                       \DeclareMathOperator \moduli {M}                 \newcommand \odd  {\ensuremath{\moduli_C(2,\cL)}}

\DeclareMathOperator\derived{\fD}               \newcommand\bounded{\ensuremath{\mathrm{b}}}    \newcommand\dbcoh[1]{\derived^\bounded(#1)}

\mathchardef\mhyphen="2D
\newcommand\dash{\nobreakdash-\hspace{0pt}}

\newcommand \Bord  {\ensuremath{\mathrm{Bord}}}   \newcommand \RBord {\ensuremath{\mathrm{RBord}}}           \newcommand \dd    {\ensuremath{\mathrm{d}}}

            \DeclareMathOperator\bessel{B}           \DeclareMathOperator\besseli{I}                                       \DeclareMathOperator\CC{C}

\DeclareMathOperator\End{End}                          

\DeclareMathOperator\HH{H}

\DeclareMathOperator\Hom{Hom}                                                                   

\DeclareMathOperator\rk{rk}                        \DeclareMathOperator\Spec{Spec}                            \DeclareMathOperator\tr{tr}

\newcommand \bC {{\ensuremath{\mathbb{C}}}}

\newcommand \bF {{\ensuremath{\mathbb{F}}}}

\newcommand \bR {{\ensuremath{\mathbb{R}}}}

\newcommand \bZ {{\ensuremath{\mathbb{Z}}}}

\newcommand \fD {{\ensuremath{\mathbf{D}}}}

\newcommand \cC {\ensuremath{\mathcal{C}}}

\newcommand \cG {\ensuremath{\mathcal{G}}}
\newcommand \cH {\ensuremath{\mathcal{H}}}

\newcommand \cK {\ensuremath{\mathcal{K}}}
\newcommand \cL {\ensuremath{\mathcal{L}}}
\newcommand \cM {\ensuremath{\mathcal{M}}}

\newcommand \rT {\ensuremath{\mathrm{T}}}

\newcommand \rZ {\ensuremath{\mathrm{Z}}}

\usepackage[margin=1.2in]{geometry}

\title{Graph potentials and topological quantum field theories}
\author{Pieter Belmans \and Sergey Galkin \and Swarnava Mukhopadhay}

\begin{document}

\maketitle

\begin{abstract}
We introduce graph potentials,
which are Laurent polynomials associated to (colored) trivalent graphs.
We show that the birational type of the graph potential
only depends on the homotopy type of the colored graph,
and use this to define a topological quantum field theory.
A similar construction was recently introduced independently
by Kontsevich--Odesskii under the name of multiplicative kernels.
We end our paper by giving an efficient computational method
to compute its partition function.
This is the first paper in a series,
and we give a survey of the applications of graph potentials
in the other parts.
\end{abstract}

\tableofcontents

\section{Introduction}
In this paper we introduce \emph{graph potentials},
a collection of Laurent polynomials
associated to trivalent graphs,
that enjoys remarkable symmetries.
\Cref{subsection:intro-graph}
is a general self-contained elementary overview,
intended as an invitation for the interested reader
to play with variations on the construction.

The main result of this paper is the construction of a novel topological quantum field theory,
the \emph{graph potential field theory}.
This result encodes the hidden structure/symmetries of graph potentials,
and explains how their invariance properties are related to the geometry
of Riemann surfaces, their degenerations and Thurston's cut systems.
The partition functions are equal to high-dimensional integrals
(the inverse Laplace transforms of the \emph{periods} for a family of level sets),
and the graph potential field theory
is the first and the only known efficient method for their computation.
For more details, see \cref{subsection:results}.

This paper is one of recommended entry points to a series of our works
on conformal field theory, mirror symmetry,
and moduli spaces of rank~2 vector bundles with fixed determinant on algebraic curves of genus~$g\geq 2$.
The relationship between graph potentials
and the geometry of moduli spaces,
as well as different aspects of mirror symmetry for these objects,
is discussed in \cref{subsection:mirror-symmetry},
where we describe the big picture for graph potentials
and the follow-up works \cite{gp-sympl,gp-decomp}.
In \cref{subsection:atiyah-floer} we speculate on how graph potentials
can be used in the context of the Atiyah--Floer conjecture.

\subsection{Graph potentials}
\label{subsection:intro-graph}
We will first introduce graph potentials in an abstract way,
without reference to other topics of interest or applications.
It would be interesting to find new choices of the input in the construction.
Our choice of input is discussed in \eqref{equation:our-choice},
and its relevance is proven by the results discussed in \cref{subsection:mirror-symmetry}.

\paragraph{General construction}
Let~$\mathbf{Y}(a,b,c)$ be a function
(defined on some space equipped with a volume form,
in our case it will be~$\mathbb{C}\setminus\{0\}$)
which is symmetric in its arguments, such that
\begin{equation}
  \label{eq:IHX}
  \mathbf{Y}(a,b,s)+\mathbf{Y}(c,d,s)
  =\mathbf{Y}(a,c,t)+\mathbf{Y}(b,d,t)
  =\mathbf{Y}(a,d,u)+\mathbf{Y}(b,c,u)
\end{equation}
for some volume-preserving transformation relating the variables.

Let~$\graph=(V,E)$ be a trivalent graph,
possibly with half-edges (or leaves).
Let~$v$ be a trivalent vertex, whose edges are labelled by~$a,b,c$, as in
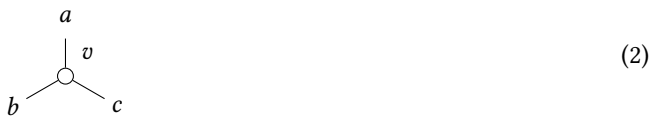
\begin{equation}
  \begin{tikzpicture}[baseline=(current bounding box.center)]
    \node[vertex] (v) at (0,0) {};
    \draw (v) + (0.3,0.3) node {$v$};

    \draw (v) edge ++ (90:0.5);
    \node (a) at (90:0.8) {$a$};

    \draw (v) edge ++ (210:0.6);
    \node (b) at (210:0.8) {$b$};

    \draw (v) edge ++ (330:0.6);
    \node (c) at (330:0.8) {$c$};
  \end{tikzpicture}
\end{equation}
We define the \emph{vertex potential}~$W_v$ as~$\mathbf{Y}(a,b,c)$.
To emphasize the dependence on~$v$ we also write~$a_v,b_v,c_v$.

Next,
the \emph{graph potential}~$W_\gamma$ is the sum of the vertex potentials associated to the (internal) vertices, i.e.,
\begin{equation}
  \label{equation:abstract-graph-potential-introduction}
  W_\graph:=\sum_{v\in V}W_v
\end{equation}
where~$W_v=W_v(a_v,b_v,c_v)$.
This is an expression in~$3(g-b_0)+2n$ arguments,
where we let~$b_0=\mathrm{b}_0(\graph)$ denote the number of connected components,
$g=\mathrm{b}_1(\graph)$ is the genus of the graph,
and~$n$ is the number of half-edges.
There are~$3(g-b_0)+n$ internal variables
and~$n$ external (or leaf) variables associated to the half-edges.

The functional equation \eqref{eq:IHX} explains how the graph potentials behave under \emph{mutation}.
Such a mutation will produce a new trivalent graph with a new graph potential associated to it,
and controlling this operation will yield important invariance properties of graph potentials,
see \cref{subsection:elementary-transformations}.
The upshot is that the different~$W_\graph$ glue to a single function~$W_{g,n}$,
and the~$W_\graph$ for different choices of~$\graph$ are restrictions to different torus charts.

\paragraph{Partition functions from graph potentials}
Applying a Fourier-like transform when the space is~$\mathbb{C}\setminus\{0\}$ we can write
\begin{equation}
  \label{equation:fourier}
  \exp(\mathbf{Y}(a,b,c))=\sum_{i,j,k}F(i,j,k)a^ib^jc^k
\end{equation}
for a function~$F(i,j,k)$. These functions satisfy the Frobenius equation
\begin{equation}
  \label{equation:frobenius-equation}
  \sum_mF(i,j,m)F(k,l,-m)=\sum_{n}F(i,k,n)F(j,l,-n).
\end{equation}
This equation hints at an important compatibility property,
culminating in \cref{theorem:alphaofalpha,theorem:tqft-introduction}.

When the space is~$\mathbb{C}\setminus\{0\}$
we can moreover consider the integration over the internal variables~$w_1,\ldots,w_{3g-3+n}$
with respect to a product of circles of sufficiently small radius around the origin,
and define the symmetric function
\begin{equation}
  Z_g(z_1,\ldots,z_n):=\iiint\exp(W_\graph(z_1,\ldots,z_n;w_1,\ldots,w_{3g-3+n})
\end{equation}
in the~$n$ external variables~$z_1,\ldots,z_n$.
This function only depends on the genus~$g$ of the graph~$\graph$,
and it satisfies the rules
\begin{equation}
  \left\{
    \begin{aligned}
      Z_{g+g'}(z_1,\ldots,z_{n+n'})
      &=\iiint Z_g(z_1,\ldots,z_n,w)Z_{g'}(z_{n+1},\ldots,z_{n+n'},w) \\
      Z_{g+1}(z_1,\ldots,z_n)
      &=\iiint Z_g(z_1,\ldots,z_n,w,w).
    \end{aligned}
  \right.
\end{equation}
This type of compatibilities,
ultimately governed by \eqref{eq:IHX},
allow us to study the partition functions associated to graph potentials.

\paragraph{Our choice of potential}
In the remainder of the paper we will focus on a specific choice of the function~$\mathbf{Y}$,
namely the Laurent polynomial
\begin{equation}
  \label{equation:our-choice}
  \mathbf{Y}(a,b,c)=\frac{a}{bc}+\frac{b}{ac}+\frac{c}{ab}+abc.
\end{equation}
The domain of the function is the torus~$(\mathbb{C}^\times)^3$.
The transformation that gives a solution to \eqref{eq:IHX} is the rational change of coordinates
\begin{equation}
  s\cdot(ab+cd)
  =t\cdot(ac+bd)
  =u\cdot(ad+bc).
\end{equation}
This transformation preserves the logarithmic rational volume forms
\begin{equation}
  \omega_a\wedge\omega_b\wedge\omega_c\wedge\omega_*\qquad *\in\{s,t,u\}
\end{equation}
where~$\omega_a=\dd\log a=\frac{\dd a}{a}$ is the Haar measure for the multiplicative group.

We can give an interpretation of this choice,
independent of the role it will play for us,
as follows.
Consider a tetrahedron inscribed in the cube~$[0,\frac{D}{2}]^3$,
whose facets are bounded by the conditions\begin{equation}
  \label{equation:inequalities}
  \left\{
    \begin{gathered}
      A\leq B+C,\qquad B\leq A+C,\qquad C\leq A+B \\
      A+B+C\leq D
    \end{gathered}
  \right.
\end{equation}
This tetrahedron can be seen as the moduli space of geodesic triangles on a sphere of radius~$R=\frac{D}{2\pi}$.
The first inequalities are the classical triangle inequalities,
and the last one is a quantum bound of the perimeter.

After the substitution~$(a,b,c,d)=(\exp A,\exp B,\exp C,\exp(-D))$
the Laurent polynomial \eqref{equation:our-choice} becomes a generating function of the inequalities \eqref{equation:inequalities},
which is invariant with respect to translations by~$2\pi\ii(A,B,C)$ for~$(A,B,C)$ satisfying
\begin{equation}
  \label{equation:integrality}
  \pm A\pm B\pm C\in\mathbb{Z}.
\end{equation}
If we let~$A,B,C$ be real (resp.~complex) numbers,
then~$a,b,c$ are real positive (resp.~complex non-zero),
and the Haar measure~$\omega_a$ in the coordinate~$a$ becomes the Lebesgue measure in the coordinate~$A$.

The definition of graph potentials for \eqref{equation:our-choice} will,
instead of trivalent graphs with half-edges,
involves colored trivalent graphs (without half-edges),
where we let a trivalent vertex~$v$ be
either uncolored~\begin{tikzpicture} \node[vertex] (A) at (0,0) {};\end{tikzpicture}
or colored \begin{tikzpicture} \node[vertex,fill] (A) at (0,0) {};\end{tikzpicture}.
We will explain this in \cref{subsection:setup-graphs}.

It would be interesting to find other functions~$\mathbf{Y}(a,b,c)$,
study their partition functions, and relate them to geometry.

\paragraph{A sneak preview of the geometry to explain our choice of potential}
As explained in \cref{subsection:mirror-symmetry}
we will relate graph potentials to the moduli space~$\odd$ of rank~2~vector bundles.
The first case to consider is for~$g=2$,
so let~$C$ be a smooth projective curve of genus~2.
In this case an explicit description of~$\odd$ due Newstead and Narasimhan--Ramanan exists \cite{MR0237500,MR0242185}:
it can be written as
\begin{equation}
  \odd\cong Q_1\cap Q_2\subset\mathbb{P}^5
\end{equation}
where~$Q_1,Q_2$ are smooth quadric hypersurfaces determined by~$C$.

These varieties degenerate to a \emph{toric} Fano threefold with six ordinary double points.
For the (smooth) intersection of quadrics in~$\mathbb{P}^5$
we have the toric degeneration given by
the (singular) intersection of quadrics with equations
\begin{equation}
  Z_0 Z_1 = Z_2 Z_3 = Z_4 Z_5 .
\end{equation}
As a toric variety it is described by the polytope
which is the convex hull of the vertices~$(\pm1,\pm1,\pm1)$.
This toric variety has terminal singularities,
and its~6 ordinary double points
admit a small resolution of singularities.

The graph potential associated to the (colored) trivalent graph of genus~$g=2$
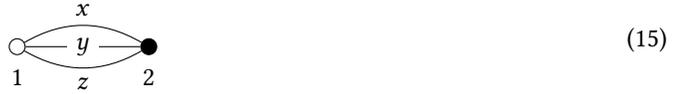
\begin{equation}
  \begin{tikzpicture}[scale=1.75,baseline=(current bounding box.center)]
    \node[vertex] (A) at (0,0) {};
    \node[vertex, fill] (B) at (1,0) {};

    \draw (A) +(0,-0.1) node [below] {$1$};
    \draw (B) +(0,-0.1) node [below] {$2$};

    \draw (A) edge [bend left]  node [above]      {$x$} (B);
    \draw (A) edge              node [fill=white] {$y$} (B);
    \draw (A) edge [bend right] node [below]      {$z$} (B);
  \end{tikzpicture}
\end{equation}
is closely related to the polytope defining the toric degeneration,
and takes on the form
\begin{equation}
  \label{equation:graph-potential-introduction}
  \widetilde{W}=xyz+\frac{x}{yz}+\frac{y}{xz}+\frac{z}{xy}+\frac{1}{xyz}+\frac{yz}{x}+\frac{xz}{y}+\frac{xy}{z}.
\end{equation}
The classical period of the Laurent polynomial~$\widetilde{W}\colon\mathbb{G}_{\mathrm{m}}^3\to\mathbb{A}^1$
agrees with the quantum period of~$\odd$,
as computed in \cite{MR3470714},
and hence the graph potential captures certain symplecto-geometric aspects of~$\odd$.
We will elaborate further on this in \cref{subsection:mirror-symmetry}.

\subsection{Overview of the results}
\label{subsection:results}
We will now discuss the results proven in this paper,
for the choice of potential in \eqref{equation:our-choice}.

\paragraph{Graph potentials under elementary transformations}
The first result we will prove
is an invariance property of graph potentials under appropriate changes of the trivalent graph.
To motivate this result,
recall that one way in which a trivalent graph of genus~$g$ arises
is as a degeneration of a smooth projective curve~$C$ of genus~$g$ into a nodal curve with only rational components.
Such a curve is also called a graph curve in Bayer--Eisenbud \cite{MR1097026},
and its dual graph gives~$\graph$.
Alternatively,
such trivalent graphs can be seen as pair of pants decompositions of~$C$,
considered as a Riemann surface.

Different degenerations,
and different pair of pants decompositions
are related to each other by elementary transformations \`a la Hatcher--Thurston,
which act on the trivalent graph.
We prove that graph potentials are suitably invariant under these elementary transformations:

\begin{alphatheorem}
  \label{theorem:alphaofalpha}
  Let~$(\graph,c)$ and~$(\graph',c')$ be two trivalent colored graphs related by elementary transformations.
  Then the graph potentials~$\widetilde{W}_{\graph,c}$ and~$\widetilde{W}_{\graph',c'}$ are identified via a rational change of coordinates.
\end{alphatheorem}

The proof of this theorem is given in \cref{subsection:elementary-transformations},
and it expresses the invariance of graph potentials under mutations.
In \cite{gp-sympl} we will further discuss the link between the mutations in \cref{theorem:alphaofalpha},
and the mutation of potentials in the context of mirror symmetry.

Associated to a Laurent polynomial we have its classical period.
The relevance of this invariant is discussed in \cref{subsection:mirror-symmetry}.
The invariance properties of graph potentials from \cref{theorem:alphaofalpha}
allow us to show in \cref{corollary:classical-periods-under-elementary-transformation}
that classical periods of graph potentials \emph{only depend} on the number of vertices and the parity of the coloring.

\paragraph{The graph potential topological quantum field theory}
Inspired by this invariance under elementary transformations,
we will define in \cref{construction:tqft} a topological quantum field theory,
and call it the \emph{graph potential field theory}.
This follows by checking the Frobenius equation \eqref{equation:frobenius-equation},
which in the context of TQFTs and their associated objects
can be rewritten as an associativity equation (for the multiplication in a Frobenius algebra)
or as the Witten--Dijkgraaf--Verlinde--Verlinde (WDVV) equation.

Its partition function~$\mathrm{Z}^{\mathrm{gp}}(t)$ (for any~$t\in\mathbb{R}$)
assigns to the pair of pants~$\Sigma_{0,3}$ with~$1+\epsilon$ anticlockwise oriented circles,
for~$\epsilon\in\{0,1\}$,
the appropriately colored vertex potential,
i.e.,
\begin{equation}
  \mathrm{Z}^{\mathrm{gp}}(t)(\Sigma_{0,3}):=\exp\left( t\widetilde{W}_{v,\epsilon}(x,y,z) \right)\in\LT^{\otimes3}.
\end{equation}
Here $\LT$ is the Hilbert space of complex valued square integrable functions on the sphere which by Foueries expansion
is same as the Hilbert space square-integrable complex series indexed by $\mathbb{Z}$.
A countable orthonormal basis is given by~$\{x^i\}_{i\in \mathbb{Z}}$
and any element of~$\LT$ can be written as $\sum_{i\in\mathbb{Z}}a_ix^{i}$ such that $\sum_{i\in\mathbb{Z}}|a_i|^2 <\infty$.
Because the potential takes values in an infinite-dimensional Hilbert space, we will need to restrict the bordism category suitably.

The properties of the graph potential can then be used to show the following theorem.
\begin{alphatheorem}
  \label{theorem:tqft-introduction}
  Let~$t\in \bC$.
  The graph potential field theory~$\mathrm{Z}^{\mathrm{gp}}(t)$
  is a two-dimensional topological quantum field theory on a suitably restricted bordism category,
  with values in Hilbert spaces.
\end{alphatheorem}
Observe that our TQFT takes values in $\LT$, and that it is not related to the usual~$(1+1)$-TQFT coming from the fusion ring of conformal blocks.

The existence of this TQFT shows why it is a powerful idea to consider moduli spaces of vector bundles on curves \emph{for all genera simultaneously}
when one is interested in mirror symmetry aspects of these Fano varieties.
In a completely different direction but with the same underlying principle,
Gonz\'alez-Prieto--Logares--Mu\~noz have constructed a topological quantum field theory
computing Hodge--Deligne polynomials for representation varieties \cite{tqft-prieto,MR4099357}.

The topological quantum field theory in \cref{theorem:tqft-introduction}
gives a powerful computational tool to compute classical periods of graph potentials.
This is discussed in \cref{section:tqft-computations}.

There is in fact an abundance of various flavours of field theories,
in this series of papers, and the works we build upon.
Indeed, in this work we have introduced a novel \emph{topological quantum field theory} to compute periods of graph potentials.
The results we build upon in \cite{gp-sympl} to construct and study toric degenerations
use the \emph{conformal field theory} given by conformal blocks \cite{MR0757857,MR0954762,MR1048605}.
On the other hand,
the study of classical and quantum periods in \cite{gp-sympl}
features both \emph{cohomological field theories}
and \emph{symplectic field theories} (see \cite{MR1291244} and \cite{MR3157146,MR2026549,MR1826267}).

Furthermore,
the cohomological and symplectic field theories in turn can be compared using a \emph{tropical field theory},
which can be seen as an intermediary between them
\cite{MR2137980,MR2404949,MR2515486,MR4068259,1902.07183}.
Finally,
at the end of this introduction
we speculate on the relationship between our work
and the construction of the (still conjectural) 4-dimensional \emph{Donaldson--Floer quantum field theory}.

\subsection{Comparison to multiplication kernels in the sense of Kontsevich--Odesskii}
\label{subsection:kontsevich-odesskii}
After the first appearance of the arXiv preprint of this article\footnote{Corresponding to Sections~2 and~3 of v1 of \cite{2009.05568v1} which was later split off to \cite{gp-tqft}.},
a similar formalism was independently developed by Kontsevich--Odesskii in \cite{MR4353308},
under the name of \emph{multiplication kernels}.
They develop the algebraic aspects more than we do,
as we mostly restrict ourselves to the specific cases
which are relevant to our geometric applications discussed
above and in \cref{section:applications}.
In what follows we will briefly explain the similarities (and differences).

In the general construction outlined in the first paragraph of \cref{subsection:intro-graph}
we consider an unspecified space with a volume form,
which in Kontsevich--Odesskii is the space~$X$.
The function~$\mathbf{Y}$ corresponds
(after the transform in \eqref{equation:fourier})
to the multiplication kernel~$K$ in op.~cit.
Thus vertex potentials correspond to multiplication kernels.
Their formalism also includes auxiliary variables
living on a priori different spaces,
which does not appear in our work.

The associativity condition from \eqref{eq:IHX}
corresponds to equation~(1.3)
(or more precisely the equation below it) in op.~cit.
The definition of a multiplication kernel of birational type,
which involves~4~variables,
makes an appearance in our \cref{subsection:tqft},
see also \eqref{equation:K4}.

Our particular choice of~$\mathbf{Y}$ in \eqref{equation:our-choice}
corresponds to their Example~2.2.
The graph potential from \eqref{equation:abstract-graph-potential-introduction}
then appears on page~25 of op.~cit.
But they only consider it for a binary rooted tree,
whereas we consider arbitrary trivalent graphs,
with or without half-edges,
along with the additional feature
of a coloring of the vertices.
Binary rooted trees correspond to trivalent graphs of genus~$g=0$
with the leaves as half-edges.

Our main motivation is to construct topological quantum field theories
whose partition functions compute period sequences
and obtain effective methods to work with them.
One of their focuses is rather to discuss algebraic structures
that arise in this context.
We briefly discuss a Frobenius algebra-like structure in \cref{subsection:tqft},
but do not develop it as they do.

It would be interesting to further compare the two formalisms,
and understand how certain other choices of multiplication kernels \`a la \cite{MR4353308}
can be used to tackle problems of the form discussed in our paper.
Observe that
our graph potentials also arise as Floer potentials
of a monotone Lagrangian torus
associated to an integrable system on the moduli space of rank-2 bundles with fixed determinant of odd degree,
as explained in \cite{gp-sympl}.
In \cite{MR4353308} similar integrable systems are an important source of multiplication kernels,
for example the integrable system associated to Hitchin systems for rank-2 bundles on~$\mathbb{P}^1$
with four (or more) punctures,
see \S3.3 and \S3.4 of op.~cit.
We leave a more detailed comparison for future work.

\paragraph{Acknowledgments}
We want to thank
Catharina Stroppel
for interesting discussions.

This collaboration started in Bonn in January--March 2018 during
the second author's visit to the
``Periods in Number Theory, Algebraic Geometry and Physics'' Trimester Program
of the Hausdorff Center for Mathematics (HIM)
and the first and third author's stay in
the Max Planck Institute for Mathematics (MPIM),
and the remaining work was done in
the Tata Institute for Fundamental Research (TIFR)
during the second author's visit in December 2019--March 2020
and the first author's visit in February 2020.
We would like to thank
HIM, MPIM and TIFR
for the very pleasant working conditions.

Several computations and experiments
were performed using
Pari/GP \cite{parigp}.

The first author was partially supported by the FWO (Research Foundation---Flanders).
The second author was partially supported by CNPq grant PQ 315747,
and Coordena\c{c}\~{a}o de Aperfei\c{c}oamento de Pessoal de N\`ivel Superior -- Brasil (CAPES) -- Finance Code 001.
The third author was partially supported by
the Department of Atomic Energy, India, under project no. 12-R\&D-TFR-5.01-0500
and also by the Science and Engineering Research Board, India (SRG/2019/000513).
 \section{Graph potentials}

\label{section:graph-potentials}
We will consider trivalent graphs of genus~$g$.
Such graphs are associated to
pair of pants decompositions of compact orientable surfaces
in the sense of Hatcher--Thurston \cite{MR0579573}.
These are collections of disjoint circles on a surface~$\Sigma$
such that their complement
is the disjoint union
of spheres with~3~holes,
also known as ``pairs of pants'' or trinions. In this section we associate a Laurent polynomial to each trivalent graph, and study how different Laurent polynomials associated to different graphs are related to each other.

The dual of the pair of pants can be encoded as follows:
\begin{equation}
  \begin{tikzpicture}[baseline=(current bounding box.center), tqft/cobordism/.style={draw}, tqft/every incoming lower boundary component/.style={draw, dashed}, tqft/every outgoing lower boundary component/.style={draw}]
    \pic[
      scale=0.5,
      tqft/pair of pants
    ];

    \node at (1.2, -0.5) {$=$};

    \node[vertex] (v) at (2,-0.5) {};
    \draw (v) edge ++ (90:0.5);
    \draw (v) edge ++ (210:0.6);
    \draw (v) edge ++ (330:0.6);
  \end{tikzpicture}
\end{equation}
and to a decomposition of~$\Sigma$ into such pairs of pants we associate the graph whose vertices are the pairs of pants, and whose edges indicate how the pairs of pants are glued together. Loops arise from cutting a genus~1~surface with 1~puncture by a circle. For every pants decomposition one thus obtains a trivalent graph.

For now we will consider trivalent graphs without making this link to geometry and topology explicit, but this will be important later.

\subsection{Definition}
\label{subsection:setup-graphs}
Let~$\graph=(V,E)$\index{G@$\graph=(V,E)$} be an undirected trivalent graph (possibly containing loops), which we will assume to be connected, and whose first Betti number is~$g$. Hence
\begin{equation}
  \begin{aligned}
    \#V&=2g-2, \\
    \#E&=3g-3.
  \end{aligned}
\end{equation}
Recall that the homology (resp.~cohomology) of a graph with coefficients in a ring~$R$ (which for us will be either~$\bZ$ or~$\bF_2$) takes~$\CC_0(\graph,R)=R^V$ and~$\CC_1(\graph,R)=R^E$ with differential given by the incidence matrix after choosing an orientation of the graph (resp.~the~$R$\dash linear dual of the incidence matrix). Hence~$\rk\HH_0(\graph,R)$ is the number of connected components, and~$\rk\HH_1(\graph,R)$ is the genus of the graph. Because we only consider the differential when~$R=\bF_2$ the choice of orientation is irrelevant for us.

We will denote
\begin{equation}
  \widetilde{N}_\graph:=\CC^1(\graph,\bZ)
\end{equation}
\index{N@$\widetilde{N}_\graph$} the free abelian group of 1-cochains on~$\graph$, and
\begin{equation}
  \widetilde{M}_\graph:=\CC_1(\graph,\bZ)
\end{equation}
\index{M@$\widetilde{M}_\graph$} the free abelian group of 1-chains.

Let~$v\in V$ be a vertex. It is adjacent to the edges~$e_{v_i},e_{v_j},e_{v_k}\in E$ (which might coincide if there is a loop). We will denote the sublattice of~$\widetilde{N}_\graph$ generated by the cochains~$x_i,x_j,x_k$ for which~$x_i(e_{v_a})=\delta_{i,v_a}$ by~$\widetilde{N}_v$\index{N@$\widetilde{N}_v$}.

Inside~$\widetilde{N}_v$ we consider the sublattice~$N_v$\index{N@$N_v$} generated by the eight cochains~$\{\pm x_i\pm x_j\pm x_k\}$. Using these we define the sublattice
\begin{equation}
  \label{equation:inclusion-of-tori}
  N_\graph\subseteq\widetilde{N}_\graph
\end{equation}
as the image of the natural morphism~$\bigoplus_{v\in V}N_v\to\widetilde{N}_\graph$.

We have the associated tori~$T_\graph^\vee:=\Spec\bC[N_\graph]$\index{T@$T_\graph^\vee$} and~$\widetilde{T}_\graph^\vee:=\Spec\bC[\widetilde{N}_\graph]$\index{T@$\widetilde{T}_\graph^\vee$}, so that~$N_\graph$ and~$\widetilde{N}_\graph$ are realized as the character lattice of~$T_{\graph}^{\vee}$ and~$\widetilde{T}_{\graph}^{\vee}$ respectively, and hence as the cocharacter lattices of~$T_\graph$ and~$\widetilde{T}_\graph$.

The inclusion~\eqref{equation:inclusion-of-tori} induces an isogeny of the tori
$\widetilde{T}_\graph^\vee\to T_\graph^\vee$, whose kernel\index{A@$A_\graph$}
\begin{equation}
  A_\graph:=\Hom(\widetilde{N}_\graph/N_\graph,\bC^\times)=(\widetilde{N}_\graph/N_\graph)^\vee,
\end{equation}
(which is isomorphic to~$\bF_2^{\oplus g}$) we wish to describe explicitly.

Consider the dual lattices~$M_\graph:= N_\graph^\vee$\index{M@$M_\graph$} and~$\widetilde{M}_\graph:=\widetilde{N}_\graph^\vee$, for which we have~$\widetilde{M}_\graph\subseteq M_\graph$, and~$A_\graph\cong M_\graph/\widetilde{M}_\graph$. The following lemma summarizes the situation.

\begin{lemma}
  We have that
  \begin{equation}
    \begin{aligned}
      M_\graph
      &=\left\{ m=\sum_{e\in E}w_ee\mid w_e\in\bR,\forall n\in N_\graph:\langle n,m\rangle\in\bZ \right\} \\
      &=\left\{ m=\sum_{e\in E}w_ee\mid w_e\in\frac{1}{2}\bZ,\forall v\in V\colon w_{e_i}+w_{e_j}+w_{e_k}\in\bZ \right\}
    \end{aligned}
  \end{equation}
  and
  \begin{equation}
    A_\graph\cong\HH_1(\graph,\bF_2).
  \end{equation}Here~$e_i, e_j, e_k$ are edges adjacent to the vertex~$v$.
\end{lemma}

In particular we also see that rational functions on the torus~$T_\graph$ are~$A_\graph$\dash invariant functions on the torus~$\widetilde{T}_\graph$.

\paragraph{Colorings}
Next we introduce colorings. We will use these to conveniently deal with a generalization of the class of trivalent graphs, where half-edges are allowed.
Using colorings we can reduce such situations (which will correspond to punctured surfaces) in the cases that we are interested in to the already considered case of trivalent graphs without half-edges.

\begin{definition}
  Let~$\graph=(V,E)$ be a trivalent graph. A \emph{coloring} is a function~$c\colon V\to\bF_2$\index{c@$c\colon V\to\bF_2$}, interpreted as an~$\bF_2$\dash valued 0-chain on~$\graph$.
\end{definition}
If~$c(v)=0$ we say that~$v$ is \emph{uncolored}, and if~$c(v)=1$ we say that~$v$ is \emph{colored}. When drawing a graph, uncolored corresponds to a circle \begin{tikzpicture} \node[vertex] (A) at (0,0) {};\end{tikzpicture} whilst colored corresponds to a disk \begin{tikzpicture} \node[vertex,fill] (A) at (0,0) {};\end{tikzpicture}. If it can be either we will indicate this by drawing \begin{tikzpicture} \node[vertex] (A) at (0,0) {}; \node[half-vertex] at (A) {};\end{tikzpicture}.

We can now introduce one of the main objects of this paper. In \cref{section:tqft-computations} we will generalize this definition to also include half-edges, but for now this definition suffices.
\begin{definition}
  \label{definition:vertex-graph-potential}
  Let~$\graph=(V,E)$ be a trivalent graph, and let~$c\colon V\to\bF_2$ be a coloring. Let~$e_1,\ldots,e_{3g-3}$ be an enumeration of the edges. We will denote~$x_i$ the coordinate variable in~$\bZ[\widetilde{N}_\graph]$ associated to~$e_i$.

  Let~$v\in V$ be a vertex, and denote~$e_i,e_j,e_k$ the three edges incident to it. Then the \emph{vertex potential} is the Laurent polynomial\index{W@$\widetilde{W}_{v,c}$}
  \begin{equation}
    \label{equation:vertex-potential}
    \widetilde{W}_{v,c}
    :=
    \hspace{-1em}\sum_{\substack{(s_i,s_j,s_k)\in\bF_2^{\oplus 3} \\ s_i+s_j+s_k=c(v)}}\hspace{-1em}x_i^{(-1)^{s_i}}\!\!\!\!x_j^{(-1)^{s_j}}x_k^{(-1)^{s_k}}
  \end{equation}
  in~$\bZ[\widetilde{N}_v]$. Observe that~$\widetilde{W}_{v,c}$ is in the image of~$\bZ[N_{v}]$ in~$\bZ[\widetilde{N}_v]$.

  Then we define the \emph{graph potential} as the Laurent polynomial\index{W@$\widetilde{W}_{\graph,c}$}
  \begin{equation}\label{equation:important-graph-potential}
    \widetilde{W}_{\graph,c}
    :=
    \sum_{v\in V}\widetilde{W}_{v,c}
  \end{equation}
  in~$\bZ[\widetilde{N}_\graph]$.~Similarly $\widetilde{W}_{\graph,c}$ is in the image of~$\bZ[N_{\graph}]$ in~$\bZ[\widetilde{N}_{\graph}]$. If~$c$ is the zero~0-chain, we will write~$\widetilde{W}_{v,0}$ and~$\widetilde{W}_{\graph,0}$.
\end{definition}

By construction we have the following result.
\begin{lemma}
  \label{lemma:descend}
  The graph potential~$\widetilde{W}_{\graph,c}$ descends to
  a regular function~$W_{\graph,c}$ on the torus $T^{\vee}_{\graph}$ (respectively the toric variety~$\widetilde{T}^{\vee}_{\graph}$)
  of the toric variety of~$X_{P_{\graph,c},M_{\graph}}$ (respectively~$X_{P_{\graph,c},\widetilde{M}_{\graph}}$).
\end{lemma}
In other words,
we have a commutative diagram
\begin{equation}
  \begin{tikzcd}
    \widetilde{T}^{\vee}_{\graph} \arrow[d] \arrow[rd, "\widetilde{W}_{\graph,c}"] \\
    T^{\vee}_{\graph} \arrow[r, swap, "W_{\graph,c}"] & \mathbb{A}^1
  \end{tikzcd}
\end{equation}

We will now make the constructions explicit.

\begin{example}
  We consider the following local picture of a trivalent vertex.
  
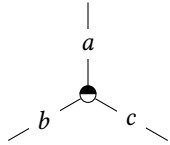
\begin{center}
    \begin{tikzpicture}[scale=1.75]
      \node[vertex] (A) at (0,0) {};
      \node[half-vertex] at (A) {};

      \draw (A) edge node [fill=white] {$a$} ++  (90:.7cm);
      \draw (A) edge node [fill=white] {$b$} ++ (210:.7cm);
      \draw (A) edge node [fill=white] {$c$} ++ (330:.7cm);
    \end{tikzpicture}
  \end{center}
  Then the vertex potentials are precisely
  \begin{equation}
    \label{equation:vertex-potentials}
    \begin{aligned}
      \widetilde{W}_{\graph,0}&=abc+\frac{a}{bc}+\frac{b}{ac}+\frac{c}{ab}, \\
      \widetilde{W}_{\graph,1}&=\frac{1}{abc}+\frac{ab}{c}+\frac{ac}{b}+\frac{bc}{a}.
    \end{aligned}
  \end{equation}
\end{example}

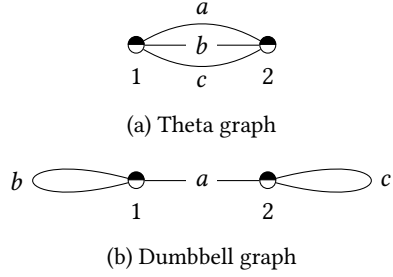
\begin{figure}[t!]
  \centering

  \begin{subfigure}[b]{\textwidth}
    \centering
    \begin{tikzpicture}[scale=1.75]
      \node[vertex] (A) at (0,0) {};
      \node[half-vertex] at (A) {};
      \node[vertex] (B) at (1,0) {};
      \node[half-vertex] at (B) {};

      \draw (A) +(0,-0.1) node [below] {$1$};
      \draw (B) +(0,-0.1) node [below] {$2$};

      \draw (A) edge [bend left]  node [above]      {$a$} (B);
      \draw (A) edge              node [fill=white] {$b$} (B);
      \draw (A) edge [bend right] node [below]      {$c$} (B);
    \end{tikzpicture}
    \caption{Theta graph}
    \label{figure:theta-graph}
  \end{subfigure}

  \begin{subfigure}[b]{\textwidth}
    \centering
    \begin{tikzpicture}[scale=1.75]
      \node[vertex] (A) at (0,0) {};
      \node[half-vertex] at (A) {};
      \node[vertex] (B) at (1,0) {};
      \node[half-vertex] at (B) {};

      \draw (A) +(0,-0.1) node [below] {$1$};
      \draw (B) +(0,-0.1) node [below] {$2$};

      \draw (A) edge [loop left]  node [left]       {$b$} (A);
      \draw (A) edge              node [fill=white] {$a$} (B);
      \draw (B) edge [loop right] node [right]      {$c$} (B);
    \end{tikzpicture}
    \caption{Dumbbell graph}
    \label{figure:dumbbell-graph}
  \end{subfigure}

  \caption{Trivalent genus 2 graphs}
  \label{figure:trivalent-g-2}
\end{figure}

There are two trivalent graphs with~2~vertices, which we will call the Theta graph and dumbbell graph respectively, and they are given in \cref{figure:theta-graph,figure:dumbbell-graph}.

\begin{example}[Theta graph]
  \label{example:theta-graph}
  We consider the Theta graph, labeled as in \cref{figure:theta-graph}. It is a trivalent graph of genus~2. We will consider two colorings for this graph, which will suffice for us by \cref{corollary:potential-only-depends-on-parity}. Let~$c$ denote the non-trivial coloring given by~$c(1)=0$, $c(2)=1$. Then we have
  \begin{align}
    \label{equation:potential-uncolored-Theta}
    \widetilde{W}_{\graph,0}
    &=2\left( abc+\frac{a}{bc}+\frac{b}{ac}+\frac{c}{ab} \right), \\
    \widetilde{W}_{\graph,c}
    &=abc+\frac{a}{bc}+\frac{b}{ac}+\frac{c}{ab} + \frac{1}{abc}+\frac{bc}{a}+\frac{ac}{b}+\frac{ab}{c}. \label{equation:theta-graph-potential-colored}
  \end{align}
\end{example}

\begin{example}[Dumbbell graph]
  \label{example:dumbbell-graph}
  Alternatively for~$g=2$ we can consider the dumbbell graph, labeled as in \cref{figure:dumbbell-graph}. We will consider two colorings for this graph, which will suffice for us by \cref{corollary:potential-only-depends-on-parity}. Let~$c$ denote the non-trivial coloring given by~$c(1)=0$, $c(2)=1$. Then we have
  \begin{align}
    \widetilde{W}_{\graph,0}
    &=ab^2+\frac{a}{b^2}+\frac{4}{a}+ac^2+\frac{a}{c^2}, \\
    \widetilde{W}_{\graph,c}
    &=ab^2+\frac{a}{b^2}+\frac{2}{a}+\frac{1}{ac^2}+2a+\frac{c^2}{a}.
  \end{align}
\end{example}

\paragraph{Independence of coloring for graph potentials}
The graph potential depends a priori on the graph~$\graph$ and the coloring $c\in\CC_0(\graph,\bF_2)$. But in fact we can identify many graph potentials as follows.

If~$e\in E$ is an edge connecting the vertices~$v$ and~$v'$ of the graph, then we have an associated~1\dash chain~$[e]\in\CC_1(\graph,\bF_2)$. Its boundary is~$\partial[e]=[v]+[v']\in\CC_0(\graph,\bF_2)$. This allows us to modify colorings~$c$ by flipping the color of~$v$ and~$v'$, preserving the parity of the coloring. Its effect on the graph potential is explained by the following lemma.

\begin{lemma}
  \label{lemma:independence-of-coloring}
  Let~$(\graph,c)$ be a graph together with a coloring~$c$.
Let~$\{x_i\}$ be a system of coordinates associated to~$\graph$.
Let~$e_k$ be an edge, and set~$c' := c+\partial[e_k]$.
Then we have the equality of graph potentials
  \begin{equation}
    \widetilde{W}_{\graph,c}=\widetilde{W}_{\graph,c'}
  \end{equation}
  after the biregular automorphism
  \begin{equation}
    x_i\mapsto
    \begin{cases}
      x_i & i\neq k \\
      x_i^{-1} & i=k
    \end{cases}
  \end{equation}
  of the algebraic torus~$\Spec\bZ[x_1^\pm,\ldots,x_{3g-3}^\pm]$.
Moreover we get a biregular automorphism of $T^{\vee}_{\graph}$
that identifies ${W}_{\graph,c}$ and $W_{\graph,c'}$.
\end{lemma}

\begin{proof}
  We only need to consider~$v$ and~$v'$, as the vertex potentials for the other vertices are not modified. It also suffices to consider~$v$. The vertex potentials are
  \begin{equation}
    \label{equation:coloring-biregular}
    \begin{aligned}
      \widetilde{W}_{v,c}
      &=\sum_{\substack{(s_i,s_j,s_k)\in\bF_2^{\oplus 3} \\ s_i+s_j+s_k=c(v)}}x_i^{(-1)^{s_i}}x_j^{(-1)^{s_j}}x_k^{(-1)^{s_k}} \\
      \widetilde{W}_{v,c'}
      &=\sum_{\substack{(s_i,s_j,s_k)\in\bF_2^{\oplus 3} \\ s_i+s_j+s_k=c(v)+1}}x_i^{(-1)^{s_i}}x_j^{(-1)^{s_j}}x_k^{(-1)^{s_k}}
    \end{aligned}
  \end{equation}
  and these sums agree after the given biregular automorphisms. More precisely, we partition~$\bF_2^{\oplus3}$ as
  \begin{equation}
    (0,0,0),(0,1,1),(1,0,1),(1,1,0)\}\sqcup\{(0,0,1),(0,1,0),(1,0,0),(1,1,1)
  \end{equation}
  and the given biregular automorphism exchanges these subsets, but then the vertex potentials are precisely identified using the given biregular automorphism.

  The second part of the proposition follows directly from the fact that each element of~$A_{\graph}$ is of order two and it acts on the variables~$x_i$'s by a character and the biregular automorphism~\eqref{equation:coloring-biregular} either fixes the variables or inverts them.
\end{proof}

Because~$\CC_1(\graph,\bF_2)$ is generated by the~1\dash chains~$[e_k]$ where~$e_k$ runs over the edges, we obtain the following corollary.
\begin{corollary}
  \label{corollary:potential-only-depends-on-parity}
  Let~$\graph$ be a trivalent graph, and~$c$ a coloring. Then the graph potential~$\widetilde{W}_{\graph,c}$ only depends on the homology class~$[c]\in\HH_0(\graph,\bF_2)$, up to biregular automorphism of the torus.
\end{corollary}

\subsection{Elementary transformations}
\label{subsection:elementary-transformations}
For a given surface of genus~$g$ there exist many (isotopy classes of) pair of pants decompositions. But they can be related via certain operations, which in \cite[Appendix]{MR0579573} are called (I), \ldots, (IV). It is remarked that (I) and (IV) in fact suffice to relate different decompositions for a surface, and that (IV) does not change the associated graph. Hence for our purposes we are only interested in the operation~(I), which we will call an \emph{elementary transformation}. Topologically this can be described as in \cref{subfigure:1a-move} on the surface and its dual graph in  \cref{subfigure:1dual-a-move}.

\begin{figure}[t!]
	\centering

	\begin{subfigure}[b]{\textwidth}
		\centering
		\begin{tikzpicture}[baseline=(current bounding box.center)]
		  \pic[
        scale=0.5,
		    rotate=90,
		    every node/.style={transform shape},
		    tqft,
		    incoming boundary components=1,
		    outgoing boundary components=2,
		    offset=-0.5,
		    cobordism edge/.style={draw},
		    every outgoing boundary component/.style={draw},
		    name=a,
		  ];
		  \pic[
        scale=0.5,
		    rotate=90,
		    every node/.style={transform shape},
		    tqft,
		    incoming boundary components=2,
		    outgoing boundary components=1,
		    anchor=outgoing boundary 1,
		    offset=0.5,
		    cobordism edge/.style={draw},
		    every incoming upper boundary component/.style={draw},
		    every incoming lower boundary component/.style={draw, dashed},
		    name=b
		  ];
		  \draw         (a-between outgoing 1 and 2) to [bend left=10]  (b-between incoming 1 and 2);
		  \draw[dashed] (a-between outgoing 1 and 2) to [bend right=10] (b-between incoming 1 and 2);
		\end{tikzpicture}
		$\longleftrightarrow$
		\begin{tikzpicture}[baseline=(current bounding box.center)]
		  \pic[
        scale=0.5,
		    rotate=90,
		    every node/.style={transform shape},
		    tqft,
		    name=a,
		    incoming boundary components=1,
		    outgoing boundary components=2,
		    offset=-0.5,
		    cobordism/.style={draw},
		    every outgoing lower boundary component/.style={draw}
		  ];
		  \pic[
        scale=0.5,
		    rotate=90,
		    every node/.style={transform shape},
		    tqft,
		    incoming boundary components=2,
		    outgoing boundary components=1,
		    offset=0.5,
		    cobordism/.style={draw},
		    anchor=outgoing boundary 1,
		    every lower boundary component/.style={draw, dashed}
		  ];
		\end{tikzpicture}
		\caption{Elementary transformation on a surface}
		\label{subfigure:1a-move}
	\end{subfigure}

	{\ }

	\begin{subfigure}[b]{\textwidth}
		\centering
		\begin{tikzpicture}[xscale=0.66,baseline=(current bounding box.center)]
		\node (i) at (0,2) [left]  {$i$};
		\node (j) at (0,0) [left]  {$j$};
		\node (k) at (3,0) [right] {$k$};
		\node (l) at (3,2) [right] {$l$};

		\draw (i) -- (1,1) -- (j) (1,1) -- (2,1) -- (k) (2,1) -- (l);
		\end{tikzpicture}
		$\longleftrightarrow$
		\begin{tikzpicture}[yscale=0.66,baseline=(current bounding box.center)]
		\node (i) at (0,3) [left]  {$i$};
		\node (j) at (0,0) [left]  {$j$};
		\node (k) at (2,0) [right] {$k$};
		\node (l) at (2,3) [right] {$l$};

		\draw (i) -- (1,2) -- (l) (1,2) -- (1,1) -- (j) (1,1) -- (k);
		\end{tikzpicture}
		\caption{Elementary transformation on the dual graph}
		\label{subfigure:1dual-a-move}
	\end{subfigure}

	\caption{Elementary transformations of trivalent graphs}
	\label{figure:1a-transformations}
\end{figure}

For later reference we summarize the discussion from \cite{MR0579573} as follows.

\begin{proposition}[Hatcher--Thurston]
  \label{proposition:hatcher-thurston}
  Elementary transformations act transitively on the set of colored trivalent graphs of genus~$g$ with~$n$ colored vertices, for~$n=0,\ldots,2g-2$.
\end{proposition}

\begin{figure}[t!]
  \centering

  \begin{subfigure}[b]{\textwidth}
    \centering
    \begin{tikzpicture}[scale=1]
      \node[vertex]       (first)  at (0,0) {};
      \node[vertex, fill] (second) at (1,0) {};
      \node[draw] (A) at (-1,1)  {};
      \node[draw] (B) at (-1,-1) {};
      \node[draw] (C) at (2,1)   {};
      \node[draw] (D) at (2,-1)  {};

      \draw (first)  +(0,-0.1) node [below] {$v_1$};
      \draw (second) +(0,-0.1) node [below] {$v_2$};

      \draw (A)     edge node [above] {$a$} (first);
      \draw (B)     edge node [below] {$b$} (first);
      \draw (first) edge node [above] {$x$} (second);
      \draw (C)     edge node [above] {$c$} (second);
      \draw (D)     edge node [below] {$d$} (second);

      \node at (3,0) {$\mapsto$};

      \node[vertex]       (first)  at (5,0.5) {};
      \node[vertex, fill] (second) at (5,-0.5) {};
      \node[draw] (A) at (4,1.5)  {};
      \node[draw] (B) at (4,-1.5) {};
      \node[draw] (C) at (6,1.5)  {};
      \node[draw] (D) at (6,-1.5) {};

      \draw (first)  node [left] {$v_1$};
      \draw (second) node [left] {$v_2$};

      \draw (A)     edge node [above] {$a$} (first);
      \draw (B)     edge node [below] {$b$} (second);
      \draw (first) edge node [right] {$x$} (second);
      \draw (C)     edge node [above] {$c$} (first);
      \draw (D)     edge node [below] {$d$} (second);
    \end{tikzpicture}
    \caption{Elementary transformation with colors}
    \label{subfigure:elementary-transformation-colors}
  \end{subfigure}

  {\ }

  \begin{subfigure}[b]{\textwidth}
    \centering
    \begin{tikzpicture}[scale=1]
      \node[vertex] (first)  at (0,0) {};
      \node[vertex] (second) at (1,0) {};
      \node[draw] (A) at (-1,1)  {};
      \node[draw] (B) at (-1,-1) {};
      \node[draw] (C) at (2,1)   {};
      \node[draw] (D) at (2,-1)  {};

      \draw (first)  +(0,-0.1) node [below] {$v_1$};
      \draw (second) +(0,-0.1) node [below] {$v_2$};

      \draw (A)     edge node [above] {$a$} (first);
      \draw (B)     edge node [below] {$b$} (first);
      \draw (first) edge node [above] {$x$} (second);
      \draw (C)     edge node [above] {$c$} (second);
      \draw (D)     edge node [below] {$d$} (second);

      \node at (3,0) {$\mapsto$};

      \node[vertex] (first)  at (5,0.5) {};
      \node[vertex] (second) at (5,-0.5) {};
      \node[draw] (A) at (4,1.5)  {};
      \node[draw] (B) at (4,-1.5) {};
      \node[draw] (C) at (6,1.5)  {};
      \node[draw] (D) at (6,-1.5) {};

      \draw (first)  node [left] {$v_1$};
      \draw (second) node [left] {$v_2$};

      \draw (A)     edge node [above] {$a$} (first);
      \draw (B)     edge node [below] {$b$} (second);
      \draw (first) edge node [right] {$x$} (second);
      \draw (C)     edge node [above] {$c$} (first);
      \draw (D)     edge node [below] {$d$} (second);
    \end{tikzpicture}
    \caption{Elementary transformation without colors}
    \label{subfigure:elementary-transformation-without-colors}
  \end{subfigure}

  \caption{Elementary transformations of trivalent graphs}
  \label{figure:elementary-transformation}
\end{figure}

Elementary transformations will give us a family of mutations of
Laurent polynomials, as discussed in \cite{MR3007265}. We are mostly interested
in the behavior of their periods under mutations induced by operations on the
trivalent graph.

When considering the behavior of an elementary transformation of a trivalent
graph on its associated graph potential, we can always split the graph potential
as
\begin{equation}
  \widetilde{W}_{\graph,c}=\widetilde{W}_{\graph,c}^{\text{mut}}+\widetilde{W}_{\graph,c}^{\text{frozen}}
\end{equation}
where the mutated part involves the variables associated to the vertices~$v_1$ and~$v_2$ adjacent to the edge corresponding to the variable~$x$. The frozen part of the graph potential is not changed, and can be ignored.

\paragraph{Elementary transformation with colors}
Let us now describe how the graph potential changes when two edges attached to vertices of different colors come together, i.e.~we consider \cref{subfigure:elementary-transformation-colors}. We will denote by~$c$ the coloring of both~$\graph$ and~$\graph'$, under the identification of the vertices. Before the transformation we have
\begin{equation}
  \begin{aligned}
    \widetilde{W}_{\graph,c}^{\text{mut}}
    &=xcd + \frac{x}{cd} + \frac{d}{cx} + \frac{d}{cx} + \frac{1}{abx} + \frac{ab}{x} + \frac{ax}{b} + \frac{bx}{a} \\
    &=\frac{1}{x}\left( ab + \frac{1}{ab} + \frac{c}{d} + \frac{d}{c} \right) + x\left( cd + \frac{1}{cd} + \frac{a}{b} + \frac{b}{a} \right).
  \end{aligned}
\end{equation}
Denoting
\begin{equation}
  \begin{aligned}
    \mu&:=\frac{1}{abcd}(c+abd)(d+abc) \\
    \nu&:=\frac{1}{abcd}(a+bcd)(b+acd)
  \end{aligned}
\end{equation}
we can write it as
\begin{equation}
  \widetilde{W}_{\graph,c}^{\text{mut}}=\frac{\mu}{x} + \nu x.
\end{equation}

After the transformation we have
\begin{equation}
  \begin{aligned}
    \widetilde{W}_{\graph',c}^{\text{mut}}
    &=x'bd + \frac{x'}{bd} + \frac{b}{dx'} + \frac{d}{bx'} + \frac{1}{acx'} + \frac{ac}{x'} + \frac{c}{ax'} + \frac{x'}{ac} \\
    &=\frac{1}{x'}\left( ac + \frac{1}{ac} + \frac{b}{d} + \frac{d}{b} \right) + x'\left( bd + \frac{1}{bd} + \frac{a}{c} + \frac{c}{a} \right).
  \end{aligned}
\end{equation}
Denoting
\begin{equation}
  \begin{aligned}
    \mu'&:=\frac{1}{abcd}(b+acd)(d+acd) \\
    \nu'&:=\frac{1}{abcd}(a+bcd)(c+abd)
  \end{aligned}
\end{equation}
we can write it as
\begin{equation}
  \widetilde{W}_{\graph',c}^{\text{mut}}=\frac{\mu'}{x'} + \nu' x'.
\end{equation}

\begin{lemma}
  \label{lemma:munu-mu'nu'}
  Let~$\mu,\nu,\mu',\nu'$ be Laurent polynomials such that~$\mu\nu=\mu'\nu'$. Then the Laurent polynomials~$\frac{\mu}{x}+\nu x$ and~$\frac{\mu'}{x'}+\nu'x'$ are identified after a rational change of coordinates.
\end{lemma}

\begin{proof}
  Setting~$z=\nu x$ we have
  \begin{equation}
    \frac{\mu}{x} + \nu x=\frac{\mu\nu}{z} + z.
  \end{equation}
  On the other hand we can do the rational change of coordinates~$z=\frac{\mu'}{x'}$ to get
  \begin{equation}
    \frac{\mu}{x}+\nu x=\frac{\mu'}{x'}+\frac{\mu\nu x'}{\mu'}.
  \end{equation}
  But by assumption we have~$\frac{\mu\nu}{\mu'}=\nu'$, hence have identified~$\frac{\mu}{x}+\nu x$ and~$\frac{\mu'}{x'}+\nu'x'$.
\end{proof}

\begin{theorem}
  \label{theorem:elementary-transformation-colors}
  Let~$\graph$ and~$\graph'$ be trivalent graphs related via an elementary transformation at an edge with \emph{different} colors. Then
  \begin{enumerate}
    \item the graph potentials~$\widetilde{W}_{\graph,c}$ and~$\widetilde{W}_{\graph',c}$ are identified after a rational change of coordinates;
    \item the rational transformation is invariant under the action of~$A_\graph$ and~$A_{\graph'}$ and hence identifies~$W_{\graph,c}$ and~$W_{\graph',c'}$.
  \end{enumerate}
\end{theorem}

\begin{proof}
  The first point follows from \cref{lemma:munu-mu'nu'}, as we have
  \begin{equation}
    \mu\nu=\mu'\nu'\frac{1}{(abcd)^2}(a+bcd)(b+acd)(c+abd)(d+abc).
  \end{equation}
  Hence the rational transformation we are using is given by~$x'=\frac{\mu'}{\nu x}$.

  To prove the second point, observe that the group~$A_\graph$ acts on the variables~$a,b,c,d,x$ via characters~$\sigma_a,\sigma_b,\sigma_c,\sigma_d,\sigma_x$, satisfying the relations
  \begin{equation}
    \label{equation:A-Gamma-relations}
    \sigma_a\sigma_b=\sigma_x=\sigma_c\sigma_d.
  \end{equation}
  Likewise~$A_{\graph'}$ acts on the variables~$a,b,c,d,x'$ via characters~$\sigma_a',\sigma_b',\sigma_c',\sigma_d',\sigma_{x'}'$, satisfying the relations
  \begin{equation}
    \label{equation:A-Gamma'-relations}
    \sigma_a'\sigma_c'=\sigma_{x'}'=\sigma_b'\sigma_d'.
  \end{equation}

  To compute the action of~$A_\graph$ and~$A_{\graph'}$ we write the rational change of coordinates more explicitly as
  \begin{equation}
    x'=\frac{(b+acd)(d+abc)}{x(a+bcd)(b+acd)}
  \end{equation}
  where we have removed the factors~$abcd$ (one coming from~$\graph$, the other from~$\graph'$) as~$\sigma_a\sigma_b\sigma_c\sigma_d$ and~$\sigma_a'\sigma_b'\sigma_c'\sigma_d'$ are trivial so they do in fact cancel. The group~$A_{\graph'}$ acts on the numerator, the group~$A_\graph$ on the denominator, and we need to check that for every~$u\in A_\graph\cong A_{\graph'}$ the action on the left- and right-hand side agrees.

  On the left-hand side we have that
  \begin{equation}
    u\cdot x'=\sigma_{x'}'(u)x'
  \end{equation}
  whilst on the right-hand side we have that
  \begin{equation}
    \begin{aligned}
      u\cdot\frac{(b+acd)(d+abc)}{x(a+bcd)(b+acd)}
      &=\frac{\left( \sigma_b'(u)b+\sigma_a'\sigma_c'\sigma_d'(u)acd \right) \left( \sigma_d'(u)d+\sigma_a'\sigma_b'\sigma_c'(u)abc \right)}{\sigma_x(u)x\left( \sigma_a(a)+\sigma_b\sigma_c\sigma_d(u)bcd \right) \left( \sigma_b(b)+\sigma_a\sigma_c\sigma_d(u)acd \right)} \\
      &=\frac{\left( \sigma_b'(u)b+\sigma_b'(u)acd \right) \left( \sigma_d'(u)d+\sigma_d'(u)abc \right)}{\sigma_x(u)x\left( \sigma_a(a)+\sigma_a(u)bcd \right) \left( \sigma_b(b)+\sigma_b(u)acd \right)} \\
      &=\frac{\sigma_b'\sigma_d'}{\sigma_x\sigma_a\sigma_b}(u)\frac{\mu'}{\nu x}.
    \end{aligned}
  \end{equation}
  But by the relations in~\eqref{equation:A-Gamma-relations} and~\eqref{equation:A-Gamma'-relations} we have the necessary equality of characters.
\end{proof}

\paragraph{Elementary transformations without colors}
Let us now describe how the graph potential changes when two edges attached to vertices of the same colors (and hence it is enough to assume that they have no coloring) come together, i.e.~we consider \cref{subfigure:elementary-transformation-without-colors}. The proof is similar to the previous case, so not all details are given. Before the transformation we have
\begin{equation}
  \begin{aligned}
    \widetilde{W}_{\graph}^{\text{mut}}
    &=xcd + \frac{x}{cd} + \frac{c}{dx} + \frac{d}{cx} + abx + \frac{a}{bx} + \frac{x}{ab} + \frac{b}{ax} \\
    &=\frac{1}{x}\left( \frac{a}{b} + \frac{b}{a} + \frac{c}{d} + \frac{d}{c} \right) + x\left( cd + \frac{1}{cd} + \frac{1}{ab} + ab \right).
  \end{aligned}
\end{equation}
Denoting
\begin{equation}
  \begin{aligned}
    \mu&:=\frac{1}{abcd}(ad+bc)(ac+bd) \\
    \nu&:=\frac{1}{abcd}(1+abcd)(cd+ab)
  \end{aligned}
\end{equation}
we can write it as
\begin{equation}
  \widetilde{W}_{\graph}^{\text{mut}}=\frac{\mu}{x}+\nu x.
\end{equation}

After the transformation we have
\begin{equation}
  \begin{aligned}
    \widetilde{W}_{\graph}^{\text{mut}}
    &=x'bd + \frac{x'}{bd} + \frac{b}{dx'} + \frac{d}{bx'} + acx' + \frac{a}{cx'} + \frac{c}{ax'} + \frac{x'}{ac} \\
    &=\frac{1}{x'}\left( \frac{b}{d} + \frac{d}{b} + \frac{a}{c} + \frac{c}{a} \right) + x\left( bd + \frac{1}{bd} + ac + \frac{1}{ac} \right).
  \end{aligned}
\end{equation}
Denoting
\begin{equation}
  \begin{aligned}
    \mu'&:=\frac{1}{abcd}(ab+cd)(ad+bc) \\
    \nu'&:=\frac{1}{abcd}(1+abcd)(ac+bd)
  \end{aligned}
\end{equation}
we can write it as
\begin{equation}
  \widetilde{W}_{\graph}^{\text{mut}}=\frac{\mu'}{x'} + \nu'x'.
\end{equation}

We obtain the following analogue of \cref{theorem:elementary-transformation-colors}, where we use that
\begin{equation}
  \mu\nu=\mu'\nu'=\frac{1}{(abcd)^2}(1+abcd)(ac+bd)(ab+cd)(ad+bc).
\end{equation}

\begin{theorem}
  \label{theorem:elementary-transformation-without-colors}
  Let~$\graph$ and~$\graph'$ be trivalent graphs related via an elementary transformation at an edge with \emph{the same} colors. Then
  \begin{enumerate}
    \item the graph potentials~$\widetilde{W}_{\graph,c}$ and~$\widetilde{W}_{\graph',c}$ are identified after a rational change of coordinates;
    \item the rational transformation is invariant under the action of~$A_\graph$ and~$A_{\graph'}$ and hence identifies~$W_{\graph,c}$ and~$W_{\graph',c}$.
  \end{enumerate}
\end{theorem}

As a direct application of \cref{theorem:elementary-transformation-colors} and \cref{theorem:elementary-transformation-without-colors}, we obtain the following corollary by \emph{tropicalizing} the rational change of coordinates.
\begin{corollary}
	Let~$(\graph,c)$ and~$(\graph',c')$ be two trivalent colored graphs related by elementary transformations, then there are piecewise linear automorphisms~${T}_{\graph,c}$ of~$\bR^{3g-3}$ that maps
	\begin{equation}
		T_{\graph,c}\colon P_{\graph,c}\rightarrow P_{\graph',c'},
	\end{equation}
  where $P_{\graph,c}$ (resp.~$P_{\graph',c'}$) are the polar duals of the Newton polytope of $\widetilde{W}_{\graph,c}$ (resp.~$\widetilde{W}_{\graph',c'}$).
\end{corollary}

\paragraph{Elementary transformations in terms of edge contractions}
More conceptually we can describe these elementary transformations in terms of edge contractions and potentials for quadrivalent vertices. 

\begin{figure}[t!]
  \centering
  \begin{tikzpicture}
    \node[vertex] (first)  at (-2,0) {};
    \node[vertex] (second) at (-1,0) {};
    \node[half-vertex] at (first)  {};
    \node[half-vertex] at (second) {};
    \node[draw] (A) at (-3,1)  {};
    \node[draw] (B) at (-3,-1) {};
    \node[draw] (C) at (0,1)   {};
    \node[draw] (D) at (0,-1)  {};

    \draw (first)  +(0,-0.1) node [below] {$v_1$};
    \draw (second) +(0,-0.1) node [below] {$v_2$};

    \draw (A)     edge node [above] {$a$} (first);
    \draw (B)     edge node [below] {$b$} (first);
    \draw (first) edge node [above] {$x$} (second);
    \draw (C)     edge node [above] {$c$} (second);
    \draw (D)     edge node [below] {$d$} (second);

    \node at (1,0) {$\mapsto$};

    \node[vertex] (first) at (3,0) {};
    \node[half-vertex] at (first) {};
    \node[draw] (A) at (2,1)  {};
    \node[draw] (B) at (2,-1) {};
    \node[draw] (C) at (4,1)  {};
    \node[draw] (D) at (4,-1) {};

    \draw (first) +(0,-0.1) node [below] {$v$};

    \draw (A) edge node [above] {$a$} (first);
    \draw (B) edge node [below] {$b$} (first);
    \draw (C) edge node [above] {$c$} (first);
    \draw (D) edge node [below] {$d$} (first);
  \end{tikzpicture}
  \caption{Edge contraction for an elementary transformation}
  \label{figure:edge-contraction}
\end{figure}
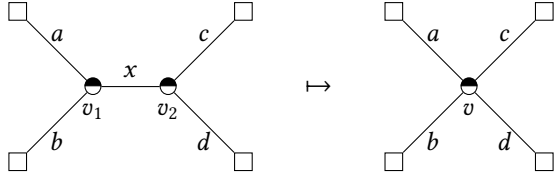

We will consider \cref{figure:elementary-transformation}, and contract the edge~$x$ between the vertices~$v_1$ and~$v_2$ into a vertex~$v$ as in \cref{figure:edge-contraction}. Associated to this vertex we add the variable~$z^{\pm}$ to the coordinate ring of the torus.

Let~$v\in V$ be the quadrivalent vertex, and denote~$a,b,c,d$ the four edges incident to it. Then the \emph{quadrivalent vertex potential} is the Laurent polynomial
\begin{equation}
  \widetilde{W}_{v,c}:=
  \begin{cases}
    \displaystyle\frac{(ab+cd)(ad+bc)(ac+bd)(1+abcd)}{(abcd)^2}\frac{1}{z}+z & c(v)=0 \\
    \displaystyle\frac{(a+bcd)(b+acd)(c+abd)(d+abc)}{(abcd)^2}\frac{1}{z}+z & c(v)=1.
  \end{cases}
\end{equation}
Here the induced coloring~$c(v)$ is defined as~$c(v_1)+c(v_2)$.

Then the mutation of the graph potential can be described as the transformation
\begin{equation}
  \widetilde{W}_{v_1,c}+\widetilde{W}_{v_2,c}\mapsto
  \frac{(ab+cd)(ad+bc)(ac+bd)(1+abcd)}{(abcd)^2}\frac{1}{z} + z
\end{equation}
if~$c(v_1)+c(v_2)=c(v)=0$, and
\begin{equation}
  \widetilde{W}_{v_1,c}+\widetilde{W}_{v_2,c}\mapsto
  \frac{(a+bcd)(b+acd)(c+abd)(d+abc)}{(abcd)^2}\frac{1}{z} + z
\end{equation}
if~$c(v_1)+c(v_2)=c(v)=1$.

\subsection{Periods of Laurent polynomials}
\label{subsection:periods-laurent}
Consider an~$n$-cycle~$\{|x_1|=\ldots=|x_n|=1\}$ in the torus~$(\bC^{\times})^{n}$, along with the normalized volume form given by~$\left( \frac{1}{2\pi\ii}  \right)^n\frac{\dd x_1}{x_1}\dots\frac{\dd x_n}{x_n}$.
\begin{definition}
  Let~$W\in\bC[x_1^{\pm},\ldots,x_n^{\pm}]$ be a Laurent polynomial and denote by~$[W]_0$ its constant term. The \emph{period}~$\pi_W(t)$\index{p@$\pi_W(t)$} of~$W$ is defined as
  \begin{equation}
    \pi_W(t):=\left( \frac{1}{2\pi\ii} \right)^n \int \dots \int_{|x_1|=\ldots=|x_n|=1} \frac{1}{1-tW} \frac{\dd x_1}{x_1}\cdots\frac{\dd x_n}{x_n},
  \end{equation}
  which can be identified with
  \begin{equation}
    \label{equation:periodpowerseries}
    \pi_W(t)=\sum_{k\geq 0}[W^k]_0t^k
  \end{equation}
  for~$|t|<1/\max\{|W(x)\mid |x|=1\}$.
\end{definition}
We will often denote the constant term~$[W^k]_0$ of the~$k$th power of~$W$ as~$\pi_k$\index{p@$\pi_k$}.

The \emph{inverse Laplace transform}~$\widecheck{\pi}_W(t)$ of~$\pi_W(t)$ is
\begin{equation}
  \widecheck{\pi}_W(t)
  :=
  \left( \frac{1}{2\pi\ii} \right)^n
  \idotsint_{|x_1|=\ldots=|x_n|=1}
  \exp(tW)\frac{\dd x_1}{x_1}\cdots\frac{\dd x_n}{x_n},
\end{equation}
which now converges absolutely and locally uniformly for all $t$.
Its everywhere convergent Taylor series expansion is given by
\begin{equation}
  \label{equation:invperiodpowerseries}
  \widecheck{\pi}_W(t)=\sum_{n\geq 0}\frac{[W^k]_0}{k!}t^k.
\end{equation}

Finally we remark that in \cref{section:tqft}, we will express~$\pi_{W}(t)$ as the trace of a trace-class operator in the Hilbert space~$\LT$.

In this way we have associated an integer sequence to a graph potential.
Let us discuss the easiest example of this, where~$g=2$.

\begin{example}[Genus two]
  \label{example:genus-two-classical-periods}
The periods for the two colorings of the Theta graph as discussed in \cref{example:theta-graph} are
  \begin{equation}
    \begin{aligned}
      \pi_{\widetilde{W}_{\graph,0}}(t)
      &=1+384t^{4}+645120t^{8}+1513881600t^{12}+\ldots \\
      \pi_{\widetilde{W}_{\graph,c}}(t)
      &=1+8t^{2}+216t^{4}+8000t^{6}+343000t^{8}+16003008t^{10}+788889024t^{12}+\ldots
    \end{aligned}
  \end{equation}
  These are in fact the same as for the dumbbell graph, which will follow from \cref{corollary:classical-periods-under-elementary-transformation}.
\end{example}

The goal of this section is to prove that the periods of the graph potential only depend on the genus of the graph and the parity of the coloring. This will follow from the following lemma, which is a mild generalization of \cite[Lemma~1]{MR3007265}, whose proof we include for completeness' sake.

\begin{lemma}
  \label{lemma:classical-periods-under-mutation}
  Let~$\varphi$ be an automorphism of the field~$\bC(x_1,\ldots,x_n)$
  corresponding to a rational transformation
  of the torus~$(\bC^\times)^n=\Spec\bC[x_1^\pm,\ldots,x_n^\pm]$,
  given by a collection of $n$
  rational functions~$(p_1/q_1,\ldots,p_n/q_n)$.
  Assume that for some~$\alpha\in\bC^\times$
  \begin{equation}
    \varphi^*\omega=\alpha\omega,
  \end{equation}
  where~$\omega=\operatorname{dlog}x_1\wedge\ldots\wedge\operatorname{dlog}x_n$ is the volume form.

  Let~$W,W'$ be Laurent polynomials in~$\bC[x_1^\pm,\ldots,x_n^\pm]$.
  If~$\varphi^*W=W'$, then~$\pi_W(t)=\pi_{W'}(t)$.
\end{lemma}

\begin{proof}
  Let~$Z_f$ be the vanishing locus of~$f=\prod_{i=1}^np_iq_i$. Consider the morphism
  \begin{equation}
    \operatorname{Log}\colon(\bC^\times)^n\to\bR^n:(z_1,\ldots,z_n)\mapsto(\log|z_1|,\ldots,\log|z_n|).
  \end{equation}
  By \cite[Corollary~6.1.8]{MR2394437} the image~$A:=\operatorname{Log}(Z_f)$, called the \emph{amoeba}, is a proper subset of~$\bR^n$,
such that~$\bR^n\setminus A$
is a disjoint union of convex sets.
There exists an element~$r\in\bR_{>0}^n$ such that $\varphi$ is regular in all points of the torus
  \begin{equation}
    T_r:=\{(z_1,\ldots,z_n)\in(\bC^\times)^n\mid \forall i=1,\ldots,n:|z_i|=r_i\} \subseteq U.
  \end{equation}
  Here~$U:=(\bC^\times)^n \backslash Z_f$.
The cycles~$T_r$,$T_{r'}$ are homologous in~$(\bC^\times)^n$ for any~$r,r'\in\bR_{>0}^n$,
hence they define the same homology class
in~$\HH_n((\bC^\times)^n,\bZ)\cong\bZ\gamma$, where~$\gamma=[T_r]$.

  This implies that $[\varphi(T_r)]=k\gamma\in\HH_n(T,\bZ)$ for some~$k\in\bZ$. But we also see that
  \begin{equation}
    \alpha
    =\alpha\int_{T_r}\omega
    =\alpha\int_{T_r}i_U^*(\omega_U)
    =\int_{T_r}\varphi_U^*(\omega_U)
    =\int_{\varphi(T_r)}\omega
    =\int_{kT_{r'}}\omega
    =k
  \end{equation}
  where~$i_U\colon U\to(\bC^\times)^n$
and~$\varphi_U$ is the composition
of the rational map~$\varphi$ with~$i_U$,
which is everywhere defined.
Hence~$\alpha=k$ is an non-zero integer.

  Now for any~$r\in\bR_{>0}^n$ such that~$T_r\subseteq U$ the equality
  \begin{equation}
      \alpha \int_{T_r}\frac{\omega}{1-tW'}
    = \int_{T_r}\varphi^*\left(\frac{\omega}{1-tW}\right)
    = \int_{\varphi(T_r)}\frac{\omega}{1-tW}
    =  k \int_{T_r} \frac{\omega}{1-tW}
  \end{equation}
  and non-vanishing of $\alpha=k$ implies the equality of periods.
\end{proof}

Hence we obtain the following corollary.

\begin{corollary}
  \label{corollary:classical-periods-under-elementary-transformation}
  Let~$(\graph,c)$ and~$(\graph',c')$ be related via elementary transformations or change of colors by the boundary of a~$1$-chain~$\CC_1(\graph,\bF_2)$. Then their periods agree, i.e.
  \begin{equation}
    \label{equation:equalityofperiods}
    \pi_{\widetilde{W}_{\graph,c}}(t)= \pi_{\widetilde{W}_{\graph',c'}}(t).
  \end{equation}
  and moreover are independent of the choice of ambient lattice, i.e.
  \begin{equation}
    \pi_{\widetilde{W}_{\graph,c}}(t)= \pi_{{W}_{\graph,c}}(t),\\
  \end{equation}
\end{corollary}

This explains why in \cref{example:genus-two-classical-periods} we could claim that the periods for the two distinct genus two graphs agree. And combining the corollary with \cref{proposition:hatcher-thurston} we have that the periods only depend on the genus and the parity of the coloring.

 \section{Topological quantum field theories from graph potentials}
\label{section:tqft}
The invariance under elementary transformations from \cref{subsection:elementary-transformations}
can be used to define a two-dimensional topological quantum field theory (or 2d~\tqft for short),
which will give an efficient computational tool to compute period sequences.
Let us quickly recall what 2d~{\tqft}s are,
using the functorial description from \textcite{MR1001453}.
For more information one is referred to, e.g., \cite{MR2037238},

A \emph{two-dimensional topological quantum field theory}
is a symmetric monoidal functor
\begin{equation}
  Z\colon(\Bord_2,\sqcup)\to(\cC,\otimes)
\end{equation}
from the symmetric monoidal category~$\Bord_2$ \index{B@$\Bord_2$} of 2-bordisms
to a symmetric monoidal category~$(\cC,\otimes)$.
An object in~$\Bord_2$ is an oriented closed curve
(a closed topological manifold of real dimension one),
i.e.~a (possibly empty)~disjoint union of copies of~$\sphere^1$,
and a morphism (or bordism) from~$E$ to~$F$
is an equivalence class of oriented compact surfaces~$M$
together with identifications of the boundary~$\partial M$ to~$E$ and~$F$.
The equivalence relation identifies bordisms via orientation-preserving diffeomorphisms,
keeping only the topological information.
We will write~$\sphere^1_{0}$ (respectively~$\sphere^1_{1}$)
for the circle with anticlockwise (respectively clockwise) orientation.

\subsection{Restricted TQFT}
For our purposes, we will take as our target category the category of complex Hilbert spaces (not necessarily finite-dimensional), with the space of bounded operators between Hilbert spaces as morphisms. This category is a monoidal category under the {\em Hilbertian tensor product} but it is not rigid. As rigidity implies the existence of the trace function for all endomorphisms which forces the Hilbert spaces to be finite-dimensional.
This also means that we cannot consider the traces of the identity operator as in \cite[page~180]{MR1001453}. As mentioned in \cite{MR1001453}, the identity operator corresponds to a  cylinder two holes  taking trace of the identity operator corresponds to the bordism given by a torus. Hence if we have a one holed torus in our bordisms category, we can not have caps or cups. Similarly we cannot have both upper and lower handles in our bordism category.

However to every Hilbert space $\cH$, we can consider its continuous dual $\cH^*$ which is linearly anti-isomorphic to $\cH$. Moreover, this gives a natural linear isomorphism between the Hilbertian tensor product~$\cH_1\otimes \cH_2$ and the space of Hilbert--Schmidt operators~$\operatorname{HS}(\cH_1^{*}, \cH_2)$.

The category of Hilbert spaces has additional structure, by the Riesz representation theorem. Namely, to every bounded operator~$f:\cH_1\rightarrow \cH_2$ we can assign the adjoint operator~$f^{\ast}:\cH_2\rightarrow \cH_1$ of a morphism, which satisfies the following properties:
\begin{enumerate}
	\item $\operatorname{Id}^{\ast}=\operatorname{Id}$;
	\item $(f\circ g)^{\ast}=g^{\ast}\circ f^{\ast}$;
  \item $f^{\ast \ast}=f$.
\end{enumerate}
Categories with morphisms satisfying the above conditions are known as $\ast$-categories. The category~$\Bord_2$ is naturally a $\ast$-category as to every bordism we can assign the opposite bordism.

If~$M$ is an object in~$\Bord_2$ and~$\overline{M}$ is the surface $M$ with opposite orientation, and we take our target category as the category of Hilbert spaces as above, then~$Z(M)$ has to be finite-dimensional. This follows from the existence of both the evaluation map~$\Hom_{\Bord_2}(\overline{M}\sqcup M,\emptyset)$ and the coevaluation map, which implies that the vector space~$Z(M)$ has a rigid dual, see for example~\cite[Proposition~1.1.8]{MR2555928}.

Hence to allow infinite-dimensional spaces we also need to slightly modify our source category, by \emph{discarding} some bordisms in~$\Bord_2$ while the objects remain the same. The new category will be denoted by~$\RBord_2$\index{R@$\RBord_2$}. Namely we will only consider bordisms given by surfaces~$\Sigma_{g,n}$ (where~$g$ is the genus and~$n$ the number of boundary components) for which the Euler characteristic~$2-2g-n$ is strictly negative, together with cylinders (to ensure we have identity morphisms), braidings and handles considered as elements of~$\Hom_{\RBord_2}(\sphere^1\sqcup \sphere^1, \emptyset)$. This will then suffice to define invariants of closed surfaces of genus~$g\geq 2$ which is important for our applications. 

We have thus defined~$\RBord_2$ as a non-full subcategory of~$\Bord_2$. Observe that by the symmetries inherent in the definition of these TQFT's (see \cref{lemma:orientation,lemma:symmetric}), the pair of pants~$\Sigma_{0,3}$ appears as a morphism~$\sphere^1\sqcup\sphere^1\to\sphere^1$ and~$\sphere^1\to\sphere^1\sqcup\sphere^1$ with the appropriate orientations, but also as a morphism~$\sphere^1\sqcup\sphere^1\sqcup\sphere^1\to\emptyset$ and~$\emptyset\to\sphere^1\sqcup\sphere^1\sqcup\sphere^1$.
\begin{remark}
  \label{remark:oppositehandles} The category $\RBord_2$ is not a $\ast$-category. Observe that we are not allowing opposite handles (nor cups or caps) in the morphisms~$\Hom_{\RBord_2}(\emptyset,\sphere^1\sqcup \sphere^1)$ in the restricted bordisms category~$\RBord_2$. A handle combined with an opposite handle gives a torus whose corresponding assignment is the trace of the identity operator. Since our target category is the category of Hilbert spaces, the trace of identity may not be defined. This is one of the main reasons for considering the restricted bordism category~$\RBord_2$.
\end{remark}

In fact, we will describe a \emph{family} of 2d~TQFT's, parametrized by~$t\in\bC$. It is only by considering the entire family of TQFT's that we can efficiently compute period sequences.

\begin{remark}
Two-dimensional TQFT's
(with values in finite-dimensional vector spaces)
can equivalently be described using Frobenius algebras
\cite{Dijkgraaf-phd,MR1414088}.
Because we consider a more general target category,
and restricted the possible bordisms,
we do not get the usual notion of a unital Frobenius algebra,
but rather we get a Hilbertian algebra.
The algebra structure on $\LT$ with an associated Spin-structure,
in the sense of \cite{MR3349955} should be analyzed.
We leave this for future work.
\end{remark}

\subsection{The graph potential TQFT}
\label{subsection:tqft}
We fix~$t\in\bC$.
Using graph potentials we will construct
a \tqft with values in Hilbert spaces for every value of~$t$.

Later we will allow~$t$ to vary,
\emph{after} we have made the identification ensuring that at least
the partition function for surfaces without punctures
is related to the period of a graph potential,
and therefore is well-behaved when we let~$t$ vary.
This suffices for our purposes.

\paragraph{\acrfull{wdvv} equations}
We can describe a 2d~\tqft in terms of the \wdvv equations.
For this, let~$(\cH,\langle-,-\rangle)$ be a Hilbert space.
This will be the value of our 2d~\tqft for~$\sphere^1_{\epsilon}$,
and later on we will take it to be~$\LT$.

We have a natural pairing~$\langle -,- \rangle_{i,j} : \cH^{\otimes n}\to\cH^{\otimes n-2}$
for all~$i<j$ by pairing the~$i$th and~$j$th factor.

Consider an assignment~$\cM_3(t)=\cM_3(t;i,j,k)\in\cH^{\otimes 3}$\index{M@$\cM_3(t)$}, where we will use~$i,j,k$ to refer to the three tensor factors. This labeling allows us to refer to specific factors in repeated tensor products of the element~$\cM_3(t)$.

Assume that~$\cM_3(t)$ is symmetric in its factors.
Then we can define the assignment\index{M@$\cM_4(t)$}
\begin{equation}
\label{equation:K4}
\cM_4(t)=\cM_4(t;i,j,k,l)
:=\langle-,-\rangle_{m,n}\left( \cM_3(t;i,j,m)\otimes\cM_3(t;k,l,n) \right)
\in\cH^{\otimes 4}.
\end{equation}
Here~$\langle-,-\rangle_{m,n}$ refers to the natural pairing of the third and sixth factor.

\begin{figure}
  \centering
  \begin{tikzpicture}[scale=1.75]
    \node[vertex] (A) at (0,0) {};
    \node[vertex] (B) at (0.5,0) {};
    \node[vertex] (C) at (3,0) {};

    \draw (A) +(0,-0.1) node [below] {$a$};
    \draw (B) +(0,-0.05) node [below] {$b$};
    \draw (C) +(0,-0.1) node [below] {$v$};

    \draw (A) -- ($(A) + (120:0.4)$);
    \draw (A) -- ($(A) + (240:0.4)$);

    \draw (B) -- ($(B) + (300:0.4)$);
    \draw (B) -- ($(B) + (60:0.4)$);

    \draw (A) -- node [above] {$e$} (B);

    \draw (C) -- ($(C) + (90:0.4)$) node [above] {$i$};
    \draw (C) -- ($(C) + (210:0.4)$) node [left] {$j$};
    \draw (C) -- ($(C) + (330:0.4)$) node [right] {$k$};
  \end{tikzpicture}
  \caption{Labeling for internal edges and trivalent vertices}
  \label{figure:labels-internal-and-trivalent}
\end{figure}
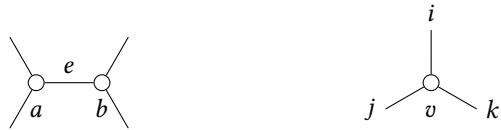

\begin{definition}
  We say that~$\cM_3(t)\in\cH^{\otimes 3}$ is a \emph{solution to the associativity equation}
  if the tensor~$\cM_4(t)\in\cH^{\otimes 4}$ is symmetric in~$i,j,k,l$.
\end{definition}
\begin{remark}
  The associativity equation encodes the associativity constraint for the multiplication in a Frobenius algebra.
  We translate it to the Frobenius equation \eqref{equation:frobenius-equation}.
  We do not equip~$\LT$ with an algebra structure as it can never be a Frobenius algebra,
  e.g., because of dimension reasons.
\end{remark}
Let~$\cM_3(t)$ be such an assignment which satisfies the associativity equation.
Then we can construct a 2d~TQFT as follows.
\begin{construction}
  \label{construction:K-Sigma-gn}
  Let~$\Sigma_{g,n}$ be an oriented surface of genus~$g$ with~$n$ punctures, such that~$2-2g-n<0$.
For every pair of pants, we will consider its dual graph. This is an oriented trivalent graph with one vertex and three half-edges. If the boundary circle is oriented anticlockwise i.e.~$\sphere^1_{0}$, then the half-edge is oriented outwards and vice versa as shown in \cref{figure:trinion-oriented}.

  \begin{figure}[h]
  	\centering
    \begin{tikzpicture}[baseline=(current bounding box.center), tqft/cobordism/.style={draw}, tqft/every incoming lower boundary component/.style={draw, dashed}, tqft/every outgoing lower boundary component/.style={draw}]
  	  \pic[scale=0.5, tqft/pair of pants, draw, name=trinion];

  	  \draw (trinion-incoming boundary 1) node [yshift=1em] {$\epsilon=0$};
  	  \draw (trinion-outgoing boundary 1) node [yshift=-1em, xshift=-.2em] {$\epsilon=1$};
  	  \draw (trinion-outgoing boundary 2) node [yshift=-1em, xshift=.2em] {$\epsilon=1$};

  	  \node at (1.7, -0.5) {$=$};

  	  \node[vertex] (v) at (3,-0.5) {};
  	  \draw[middlearrow={>}] (v) -- ++ (90:.75);
  	  \draw[middlearrow={<}] (v) -- ++ (210:.9);
  	  \draw[middlearrow={<}] (v) -- ++ (330:.9);
  	\end{tikzpicture}
    \caption{Dual graph to oriented pair of pants}
    \label{figure:trinion-oriented}
  \end{figure}
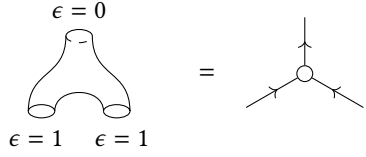

  By choosing a pair of pants decomposition of~$\Sigma_{g,n}$ we assign the dual graph~$\graph$. Observe that~$\graph$ has oriented half-edges and internal edges are unoriented. We can assign\index{M@$\cM_{\Sigma_{g,n}}$}
  \begin{equation}
    \label{equation:formal-tqft}
    \cM_{\Sigma_{g,n}}(t)
    :=
    \bigotimes_{e\in E_{\mathrm{int}}}\langle-,-\rangle_{a,b}\left( \bigotimes_{v\in V}\cM_3(t;i,j,k) \right)\in\cH^{\otimes n}
  \end{equation}
  where we use the labeling for internal edges and trivalent vertices as in \cref{figure:labels-internal-and-trivalent}, and the tensor product over the internal edges means we apply all possible pairings~$\langle-,-\rangle_{a,b}$ where the vertices~$a$ and~$b$ refer to specific factors in the tensor product~$\bigotimes_{v\in V}\cM_3(t;i,j,k)$.
\end{construction}

Using the relationship between 2d~TQFT's and solutions to the associativity equation
(see, e.g.,~\cite[\S9]{MR1492534}, where the WDVV equation is referred to as the associativity equation)
we obtain the following result.

\begin{proposition}
  \label{proposition:tqft-from-wdvv}
  Let~$\cM_3(t)\in\cH^{\otimes 3}$ be a solution to the associativity equation.
  Then \cref{construction:K-Sigma-gn} is independent of the choice of the pair of pants decomposition. \end{proposition}

\paragraph{Graph potentials as solutions to the associativity equation}
Now we will revisit the setting introduced in the introduction, in particular the functional equation in \eqref{eq:IHX}
and the ensuing discussion.
We set~$\cH$ to be~$\LT$.
By Fourier expansion, any~$f \in \LT$ can be written as
\begin{equation}
  \sum_{i\in\bZ}a_iz^i,
\end{equation}where the collection~$\{z^i\}_{i\in \bZ}$ is a complete orthonormal basis of~$\LT$ with respect to the standard pairing on~$\LT$ given by
\begin{equation}
  \label{equation:explicit-pairing}
  \begin{aligned}
    \left \langle f(z),g(z)\right\rangle_z
    &=\left\langle \sum_{i\in\bZ}a_iz^i, \sum_{j\in\bZ}b_jz^j \right\rangle_z\\
    &=\sum_{i\in\bZ}a_i\overline{b}_i\\
    &=\frac{1}{2\pi\ii}\int_{\sphere^1}f(z)\overline{g(z)}\frac{\dd z}{z}.
  \end{aligned}
\end{equation}

The associativity equation can be rephrased as follows in this setting. Consider a~$f(x_1,x_2,x_3)\in\LT^{\otimes3}$. Assume now that there exists a function~$\varphi$ such that
\begin{equation}
  f(x_1,x_2,y)f(x_3,x_4,y)=f(x_1,x_3,z)f(x_2,x_4,z)
\end{equation}
in~$\LT^{\otimes6}$, where~$y=\varphi(z,x_1,x_2,x_3,x_4)$ is a function such that the Jacobian of~$\varphi$ is the identity. In that case we get that~$\varphi_*\dlog z=\dlog y$, and hence
\begin{equation}
  \begin{aligned}
    \langle f(x_1,x_2,y)\otimes  f(x_3,x_4,y^{-1})\rangle_{y}
    &=\int_{\sphere^1}f(x_1,x_2,y)f(x_3,x_4,y)\frac{\dd y}{y}\\
    &=\int_{\sphere^1}f(x_1,x_3,z)f(x_2,x_4,z)\frac{\dd z}{z}\\
    &=\langle f(x_1,x_3,z)\otimes f(x_2,x_4,z^{-1})\rangle_{z}
  \end{aligned}
\end{equation}
in~$\LT^{\otimes 4}$ giving a solution to the associativity equation. Taking the logarithm of the function~$f$ we can interpret the multiplicative condition as an \emph{additive} condition, and this brings us in the setting of \cref{section:graph-potentials}.

In~\eqref{equation:vertex-potential} we have defined the vertex potential at a vertex~$v\in V$ of a colored graph~$(\graph,c)$. Using the variables~$x,y,z$, and writing~$\epsilon=c(v)$ we have that
\begin{equation}
  \widetilde{W}_{v,\epsilon}
  =(xyz)^{(-1)^\epsilon}+(xy^{-1}z^{-1})^{(-1)^\epsilon}+(x^{-1}yz^{-1})^{(-1)^\epsilon}+(x^{-1}y^{-1}z)^{(-1)^\epsilon}.
\end{equation}
We observe that~$\widetilde{W}_{v,\epsilon}$ is symmetric in the variables~$x,y,z$. We also have the following symmetries, aside from the symmetry in the variables.
\begin{lemma}
  \label{lemma:orientation}
  The vertex potential satisfies
  \begin{equation}
    \begin{aligned}
      \widetilde{W}_{v,0}(x,y,z)
      &=\widetilde{W}_{v,0}(x^{-1},y^{-1},z) \\
      &=\widetilde{W}_{v,1}(x^{-1},y,z) \\
      \widetilde{W}_{v,1}(x,y,z)
      &=\widetilde{W}_{v,0}(x^{-1},y,z) \\
      &=\widetilde{W}_{v,0}(x^{-1},y^{-1},z^{-1}).
    \end{aligned}
  \end{equation}
\end{lemma}

\begin{definition}
  Let~$\epsilon$ be~$0$ or~$1$, and~$t\in\bC$. We define
  \begin{equation}
    \label{equation:f-epsilon}
    \begin{aligned}
      f_{\epsilon}(x,y,z;t)
      &:=\exp\left( t\widetilde{W}_{v,\epsilon}(x,y,z) \right) \\
      &=\sum_{a,b,c,d\geq 0}\frac{t^{a+b+c+d}}{a!b!c!d!}x^{(-1)^\epsilon\left( (a+b)-(c+d) \right)}y^{(-1)^\epsilon\left( (a+c)-(b+d) \right)}z^{(-1)^\epsilon\left( (a+d)-(b+c) \right)}.
    \end{aligned}
  \end{equation}
\end{definition}

The following lemma follows directly. \begin{lemma}
  \label{lemma:symmetric}
  For all~$t\in\bC$ we have that~$f_\epsilon(x,y,z;t)\in\LT^{\otimes3}$. It is moreover symmetric in~$x,y,z$.
\end{lemma}

Finally, as observed above, we have translated between an additive and a multiplicative form of the associativity equation, and hence by \cref{subsection:elementary-transformations} we obtain the following
\begin{proposition}
  For all~$t\in\bC$ we have that~$f_\epsilon(x,y,z;t)$ satisfies the associativity equation.
\end{proposition}

We now come to the essential construction.
\begin{construction}
  \label{construction:tqft}
  Let~$t\in\bC$. The \emph{graph potential field theory}~$\rZ^{\mathrm{gp}}(t)$\index{Z@$\rZ^{\mathrm{gp}}(t)$} for~$t$ is defined as follows.

  Let~$\Sigma_{0,3;\epsilon_1,\epsilon_{2},\epsilon_{3}}$ be a pair of pants, and $\epsilon_1,\epsilon_2, \epsilon_3$ denote the orientation of the three boundary circles. Let~$\graph$ be the trivalent graph on one vertex and three oriented half-edges which is dual to the pair of pants. Then we define the partition function as
  \begin{equation}\label{equation:paritymatchpants}
    \rZ^{\mathrm{gp}}(t)(\Sigma_{0,3;\epsilon_1,\epsilon_{2},\epsilon_{3}}):= \exp\left( t\widetilde{W}_{v,0}(x^{(-1)^{\epsilon_1}},y^{(-1)^{\epsilon_2}},z^{(-1)^{\epsilon_3}}) \right)
  \end{equation}
  in~$\LT^{\otimes3}$. Here~$v$ is the unique vertex of the dual graph~$\graph$ of~$\Sigma_{0,3;\epsilon_1,\epsilon_{2},\epsilon_{3}}$. We assign~$x,y,z$ (resp.~$x^{-1}$, $y^{-1}$, $z^{-1}$) as coordinate variables corresponding to the half-edges that are oriented outwards (resp.~inwards) as shown in \cref{figure:variableattachement}.

  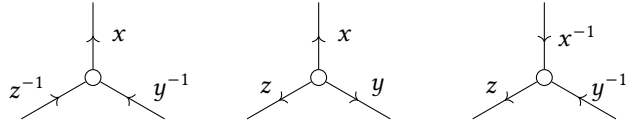
\begin{figure}[h]
  	\centering
  	\begin{tikzpicture}
  	\node[vertex] (v) at (0,0) {};
  	\draw[middlearrow={>}] (v) -- node [xshift=1em] {$x$} ++ (90:1);
  	\draw[middlearrow={<}] (v) -- node [xshift=-1em, yshift=.5em] {$z^{-1}$} ++ (210:1.1);
  	\draw[middlearrow={<}] (v) -- node [xshift=1.5em, yshift=.5em] {$y^{-1}$} ++ (330:1.1);

  	\node[vertex] (w) at (3,0) {};
  	\draw[middlearrow={>}] (w) -- node [xshift=1em] {$x$} ++ (90:1);
  	\draw[middlearrow={>}] (w) -- node [xshift=-.5em, yshift=.5em] {$z$} ++ (210:1.1);
  	\draw[middlearrow={>}] (w) -- node [xshift=.75em, yshift=.5em] {$y$} ++ (330:1.1);

  	\node[vertex] (z) at (6,0) {};
    \draw[middlearrow={<}] (z) -- node [xshift=1.25em] {$x^{-1}$} ++ (90:1);
  	\draw[middlearrow={>}] (z) -- node [xshift=-.5em, yshift=.5em] {$z$} ++ (210:1.1);
  	\draw[middlearrow={<}] (z) -- node [xshift=1em, yshift=.5em] {$y^{-1}$} ++ (330:1.1);
  	\end{tikzpicture}
  	\caption{Variable attachment to oriented graphs}
  	\label{figure:variableattachement}
  \end{figure}

    The assignment to the oriented pair of pants is well-defined by virtue of \cref{lemma:orientation}. Moreover  \cref{lemma:orientation} also implies that
  \begin{equation}\label{equation:sneakylemma}
   \rZ^{\mathrm{gp}}(t)(\Sigma_{0,3;\epsilon_1,\epsilon_{2},\epsilon_{3}})= \exp\left( t\widetilde{W}_{v,\epsilon}(x,y,z) \right),
  \end{equation}
where $\epsilon=\epsilon_1+\epsilon_2+\epsilon_3$.

Hence all the graphs shown in \cref{figure:variableattachement} have the same partition function~$\exp\left( t\widetilde{W}_{v,0}(x,y,z) \right)$. Similarly all the graphs in \cref{figure:variableattachementwithoopirentation} have the same partition function~$\exp\left( t\widetilde{W}_{v,1}(x,y,z) \right)$.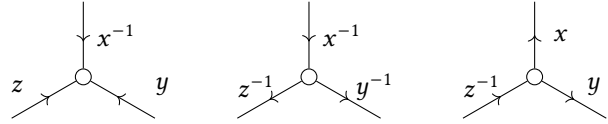
\begin{figure}[h]
	\centering
	\begin{tikzpicture}
	\node[vertex] (v) at (0,0) {};
  \draw[middlearrow={<}] (v) -- node [xshift=1.25em] {$x^{-1}$} ++ (90:1);
	\draw[middlearrow={<}] (v) -- node [xshift=-1em, yshift=.5em] {$z$} ++ (210:1.1);
	\draw[middlearrow={<}] (v) -- node [xshift=1.5em, yshift=.5em] {$y$} ++ (330:1.1);

	\node[vertex] (w) at (3,0) {};
  \draw[middlearrow={<}] (w) -- node [xshift=1.25em] {$x^{-1}$} ++ (90:1);
	\draw[middlearrow={>}] (w) -- node [xshift=-.5em, yshift=.5em] {$z^{-1}$} ++ (210:1.1);
	\draw[middlearrow={>}] (w) -- node [xshift=1em, yshift=.5em] {$y^{-1}$} ++ (330:1.1);

		\node[vertex] (z) at (6,0) {};
	\draw[middlearrow={>}] (z) -- node [xshift=1em] {$x$} ++ (90:1);
	\draw[middlearrow={<}] (z) -- node [xshift=-.5em, yshift=.5em] {$z^{-1}$} ++ (210:1.1);
	\draw[middlearrow={>}] (z) -- node [xshift=.75em, yshift=.5em] {$y$} ++ (330:1.1);
	\end{tikzpicture}
	\caption{Variable attachment to oppositely oriented graphs}
	\label{figure:variableattachementwithoopirentation}
\end{figure}

This is an additional feature of our partition function and  we color the vertex of~$\graph$ with color~$\epsilon=\epsilon_1+\epsilon_2+\epsilon_3$ to mark this feature.

  To the cylinder in~$\Hom_{\RBord_2}(\sphere^1_{\epsilon}, \sphere^1_{\epsilon+1})$ we assign the identity morphism in~$\End(\LT)$, and similarly for any braiding we just assign the identity morphism on all the factors. The handle in~$\Hom_{\RBord_2}(\sphere_{\epsilon}^1\sqcup \sphere^1_{\epsilon+1},\emptyset)$ is assigned the evaluation map which is also natural pairing on the Hilbert space $\LT$. \iffalse
   Denoting~$\langle-,-\rangle_z$ the component of the pairing in the variable~$z$, we have that
  \begin{equation}
    \begin{aligned}
    Z(\Xi)\circ Z(\Sigma_{0,3,\epsilon_1,\epsilon_2,\epsilon_3})
&=      T\circ \exp( t\widetilde{W}_{v,0}(x^{(-1)^{\epsilon}},y^{(-1)^{\epsilon}},z^{(-1)^{\epsilon}}) )\\
      &=\exp( t\widetilde{W}_{v,0}(x^{(-1)^{\epsilon+1}},y^{(-1)^{\epsilon}},z^{(-1)^{\epsilon}}) )\\
      &=Z(\Sigma_{0,3,\epsilon_1+1,\epsilon_2,\epsilon_3})\\
      \left\langle \exp\left( t\widetilde{W}_{v,\epsilon}(x,y,z) \right),\sum_{i\in\bZ}z^iw^{i} \right\rangle_{z}
      &=\exp( t\widetilde{W}_{v,\epsilon}(x,y,w)).
    \end{aligned}
  \end{equation}
  This implies that attaching a cylinder has no effect on the partition function of pants, hence it goes to the identity in~$\End(\LT^{\otimes 3})$. Similarly attaching a handle with identifying one~$\sphere^1$ from the handle and another from the pants changes the parity of the partition function of the pair of pants. In particular our assignment of the partition function to the pants is consistent with attaching handles.
\fi

  Let~$\Sigma_{g,n}$ be a connected oriented surface of genus~$g$ with~$n$ boundary components satisfying~$2-2g-n<0$. Consider a pair of pants decomposition for~$\Sigma_{g,n}$, and let~$\graph$ be the trivalent dual graph with~$n$ half-edges determined by the pair of pants decomposition. We incorporate the orientation of~$\Sigma_{g,n}$ as follows. Cut~$\graph$ along the internal edges~$E_{\mathrm{int}}$ to form a forest consisting of~$2g-2$ trivalent graphs with one vertex and three half-edges with appropriate orientations. The coloring of the vertices is determined by the parity of the orientations of the number of clockwise circles of each pair of pants. We set
  \begin{equation}
    \rZ^{\mathrm{gp}}(t)(\Sigma_{g,n}):=
    \bigotimes_{e\in E_{\mathrm{int}}}\langle-,-\rangle_{a,b}\left( \bigotimes_{v\in V}\exp\left( t\widetilde{W}_{v,\epsilon}(x,y,z) \right) \right)
  \end{equation}
  in~$\LT^{\otimes n}$. As in~\eqref{equation:formal-tqft}, we use the labeling for internal edges and trivalent vertices as in \cref{figure:labels-internal-and-trivalent}, and the tensor product over the internal edges means we apply all possible pairings~$\langle-,-\rangle_{a,b}$ where the vertices~$a$ and~$b$ refer to specific factors in the tensor product indexed by the set of vertices of~$V$.
\end{construction}

\begin{figure}[h]
	\centering
	\begin{tikzpicture}[baseline=(current bounding box.center)]
	\pic[
	scale=.5,
	rotate=90,
	every node/.style={transform shape},
	tqft,
	incoming boundary components=2,
	outgoing boundary components=0,
	offset=0.5,
	cobordism edge/.style={draw},
	name=a,
	];

	\pic[
	scale=.5,
	rotate=90,
	every node/.style={transform shape},
	tqft,
	incoming boundary components=1,
	outgoing boundary components=2,
	offset=-0.5,
	anchor=outgoing boundary 1,
	cobordism edge/.style={draw},
  every incoming lower boundary component/.style={draw, dashed},
	every incoming upper boundary component/.style={draw},
	name=b,
	];

  \draw[dashed] (0,.5) circle (12pt);

	\node at (1.5, 0.5) {$=$};

	\pic[
	xshift=10em,
	scale=.5,
	rotate=90,
	every node/.style={transform shape},
	tqft,
	incoming boundary components=2,
	outgoing boundary components=0,
	offset=0.5,
	cobordism edge/.style={draw},
	name=c,
  incoming lower boundary component 2/.style={dashed, draw},
	incoming upper boundary component 2/.style={dashed, draw},
	];

	\pic[
	xshift=10em,
	scale=.5,
	rotate=90,
	every node/.style={transform shape},
	tqft,
	incoming boundary components=1,
	outgoing boundary components=2,
	offset=-0.5,
	anchor=outgoing boundary 1,
	cobordism edge/.style={draw},
  every incoming lower boundary component/.style={draw, dashed},
	every incoming upper boundary component/.style={draw},
	name=d,
	];

	\draw (b-incoming boundary 1) node [xshift=-1em, yshift=-1em] {$\epsilon$};
	\draw (d-incoming boundary 1) node [xshift=-1em, yshift=-1em] {$\epsilon$};
	\end{tikzpicture}
	\caption{Fundamental relation for one-holed torus}
	\label{figure:cuttingandgluingalongdifferentcircles}
\end{figure}
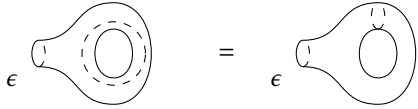

By~\eqref{equation:paritymatchpants},~\eqref{equation:sneakylemma} and \cref{lemma:orientation}, we get
\begin{equation}
  W_{v,0}(x^{(-1)^{\epsilon_1}},y^{(-1)^{\epsilon_2}},z^{(-1)^{\epsilon_3}})=W_{v,0}(x^{(-1)^{\epsilon_1}},y^{(-1)^{\epsilon_2+1}},z^{(-1)^{\epsilon_3+1}}).
\end{equation}
Equating the variables $y$ and $z$ and taking inner product in the variable $y$, we get
\begin{equation}
  \langle W_{v,0}(x^{(-1)^{\epsilon_1}},y^{(-1)^\epsilon},y^{(-1)^{\epsilon+1}})\rangle_y
  =\langle W_{v,0}(x^{(-1)^{\epsilon_1}},y^{(-1)^{\epsilon+1}},y^{(-1)^{\epsilon}})\rangle_y
\end{equation}
This equality encodes that we obtain the same partition function for one-holed torus with decompositions obtained from cutting and gluing with two different circles as shown in \cref{figure:cuttingandgluingalongdifferentcircles}. Thus the partition function $\rZ^{\operatorname{gp}}(t)(\Sigma_{1,1})$ for a one-holed torus is well-defined.

Hence by \cref{proposition:tqft-from-wdvv}, the fact that cylinders and braidings are the identity, and that the handle is assigned the natural pairing, we obtain the following
\begin{corollary}
For all~$t \in \bC$, the assignment~$\rZ^{\mathrm{gp}}(t)$ defines a two-dimensional \tqft on the restricted bordism category~$\RBord_2$.
\end{corollary}

This proves \cref{theorem:tqft-introduction}.

It would be interesting to formalize this notion of a family of infinite-dimensional TQFT's, so that one can consider all of them together. We will not develop such a formalism. Rather we will only consider the family of partition functions for closed surfaces, defined by this TQFT. This will allow us in \cref{section:tqft-computations} to obtain an efficient method to compute periods of graph potentials.

 \section{Computing periods via the graph potential TQFT}
\label{section:tqft-computations}
We can now turn our discussion to a practical method to compute the periods of graph potentials using topological quantum field theories. We also refer to \cite{MR1415320,MR2322178} for various interesting interpretations of periods, and the use of infinite-dimensional Hilbert spaces (albeit using different methods) to compute them.

The following definition extends \cref{definition:vertex-graph-potential}, where we now allow leaves (or half-edges) to be present in the graph, as the definition of the vertex potentials does not depend on whether we have half-edges or not.
\begin{definition}
  Let~$(\graph,c)$ be a colored trivalent graph of genus~$g$ with~$n$ leaves. We define the \emph{graph potential}~$\widetilde{W}_{\graph,c}$ as the sum of vertex potentials.
\end{definition}

For a graph with half-edges we let~$E_{\mathrm{int}}$ denote the set of internal edges, i.e.~we remove any half-edges from~$E$. Let us enumerate the variables associated to internal edges~$E_{\mathrm{int}}$ as~$x_1,\ldots,x_{3g-3+n}$. We introduce the following notation.
\begin{notation}
  Let~$(\graph,c)$ be a colored trivalent graph of genus~$g$ with~$n$ leaves such that~$2-2g-n<0$. Orient the half-edges of~$\graph$ such that half-edges attached to vertices~$v$ of color~$c(v)=0$ are pointing outwards while those attached to vertices of~$c(v)=1$ are pointing inwards. This orientation is consistent with the orientation coming from pair of pants decomposition.

  Denote by\index{K@$\cK_{\graph,c}(t)$}
  \begin{equation}
    \label{equation:K-with-leaves}
    \cK_{\graph,c}(t)
    :=
    \left( \frac{1}{2\pi\ii} \right)^{\#E_{\mathrm{int}}}
    \idotsint_{(\sphere^1)^{E_{\mathrm{int}}}}\exp\left( t\widetilde{W}_{\graph,c}(x_1,\ldots,x_{3g-3+2n}) \right)\frac{\dd x_1}{x_1}\cdots\frac{\dd x_{3g-3+n}}{x_{3g-3+n}}
  \end{equation}
  the corresponding element of~$\LT^{\otimes n}$.
\end{notation}

We record the following important observation as a lemma.
\begin{lemma}
  \label{lemma:laplace-transform-interpretation}
  If~$n=0$, then~\eqref{equation:K-with-leaves} reduces to
  \begin{equation}
    \cK_{\graph,c}(t)
    =\widecheck{\pi}_{\widetilde{W}_{\graph,c}}(t)
  \end{equation}
  where~$\widecheck{\pi}_{\widetilde{W}_{\graph,c}}(t)$ is the inverse Fourier--Laplace transform of the period~$\pi_{\widetilde{W}_{\graph,c}}(t)$ of the graph potential.
\end{lemma}
In other words, if~$\pi_{\widetilde{W}_{\graph,c}}(t)=\sum_{n\geq 0}\pi_nt^n$, then~$\cK_{\graph,c}(t)=\sum_{n\geq 0}p_nt^n$, where~$p_n=\pi_n/n!$.

\begin{figure}
  \centering
  \begin{tikzpicture}[scale=1.75]
    \node[vertex] (A) at (0,0) {};
    \node[vertex] (B) at (0.5,0) {};
    \node[vertex] (Ap) at (3,0) {};
    \node[vertex] (Bp) at (4,0) {};

    \draw (A) +(0,-0.1) node [below] {$a$};
    \draw (B) +(0,-0.05) node [below] {$b$};
    \draw (Ap) +(0,-0.1) node [below] {$a$};
    \draw (Bp) +(0,-0.05) node [below] {$b$};

    \draw (A)  -- ($(A) + (120:0.4)$);
    \draw (A)  -- ($(A) + (240:0.4)$);
    \draw (Ap) -- ($(Ap) + (120:0.4)$);
    \draw (Ap) -- ($(Ap) + (240:0.4)$);
    \draw (Ap) -- node [above] {$e'$} ($(Ap) + (0:0.4)$);

    \draw (B)  -- ($(B) + (300:0.4)$);
    \draw (B)  -- ($(B) + (60:0.4)$);
    \draw (Bp) -- ($(Bp) + (300:0.4)$);
    \draw (Bp) -- ($(Bp) + (60:0.4)$);
    \draw (Bp) -- node [above] {$e''$} ($(Bp) + (180:0.4)$);

    \draw (A)  -- node [above] {$e$} (B);
  \end{tikzpicture}
  \caption{Cutting an internal edge}
  \label{figure:cutting-internal-edge}
\end{figure}
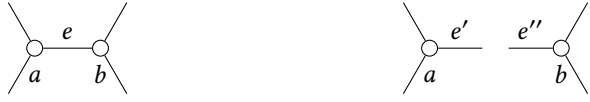

We have the following proposition, which follows from the change of variables formula for integrals. It is an important computational tool in what follows.
\begin{proposition}
  \label{proposition:cutting-a-graph}
  Let~$\graph'$ be a trivalent graph with~$n+2$ leaves. Let~$c$ be a coloring. Consider half-edges~$e'$ and~$e''$ adjacent to vertices~$a$ and~$b$. We define a new colored trivalent graph~$\graph$ (with~$n$ leaves) by replacing two leaves at the vertices~$a$ and~$b$ by the internal edge~$e$ connecting~$a$ and~$b$, as in \cref{figure:cutting-internal-edge}. Then
  \begin{equation}
    \cK_{\graph,c}(t)
    =
    \frac{1}{2\pi\ii}\int_{\sphere^1}\cK_{\graph',c}(t)|_{x_{e'}=x_{e''}=z}\frac{\dd z}{z}.
  \end{equation}
  where~$x_{e'}$ and~$x_{e''}$ are variables associated to the leaves~$e'$ and~$e''$ attached to the vertices~$a$ and~$b$.
\end{proposition}

Observe that not all trivalent colored graphs that we considered in \cref{section:graph-potentials} arise as the dual graph of a pair of pants decomposition of a orientable surface whose boundary has induced orientations. Since the category $\RBord_2$ only consists of objects of this form, we need the following results to use the TQFT partition function $\rZ^{\operatorname{gp}}(t)$ effectively to compute periods of arbitrary $(\graph,c)$.

The following proposition relates~$\cK_{\graph,c}(t)$ to the partition function of the~$\LT$-valued TQFT that we constructed.
\begin{proposition}\label{proposition:cuttingandparitymatching}
	Let~$\Sigma_{g,n}$ be a oriented surface of genus~$g$ with~$n$ boundary components and let~$(\graph,c)$ be the dual graph obtained from a pair of pants decomposition of~$\Sigma_{g,n}$. If~$3g-3+n$ is even, then
	\begin{equation}
    \cK_{\graph,c}(t)=\rZ^{\operatorname{gp}}(t)(\Sigma_{g,n}).
	\end{equation}
\end{proposition}
\begin{proof}
  \cref{proposition:cutting-a-graph} allows us to write $\cK_{\graph,c}(t)$ as the iterated integral over the pair of pants decomposition. Since we need $3g-3+n$ cuts to get to the pair of pants, the parity being even guarantees that we can use \cref{lemma:orientation} to match up the integral with the norms in $\LT$ that appears in the definition of $\rZ^{\operatorname{gp}}(t)(\Sigma_{g,n})$.
\end{proof}

\begin{remark}\label{remark:klugscheisser}
  If~$3g-3+n$ is odd, then we can compute~$\cK_{\graph,c}(t)$ by first cutting~$(\graph,c)$ along one edge to produce a new graph $(\graph',c')$, and then use \cref{proposition:cuttingandparitymatching} for $(\graph',c')$  and apply \cref{proposition:cutting-a-graph}.
\end{remark}

\paragraph{Bessel functions}
To adequately work with the partition functions of this topological quantum field theory we recall that the modified Bessel function of the second kind is defined as\index{I@$\besseli_\alpha(z)$}
\begin{equation} \label{d:besseli}
  \besseli_\alpha(z):=\sum_{m\geq 0}\frac{1}{m!\Gamma(m+\alpha+1)}\left( \frac{z}{2} \right)^{2m+\alpha}.
\end{equation}
For our purposes we are only interested in the case~$\alpha=0$, with a rescaling of the argument. We will use the following notation.
\begin{notation}
  We denote\index{B@$\bessel(z)$}
  \begin{equation} \label{d:bessel}
    \bessel(z):=\besseli_{\alpha=0}(2z)=\sum_{m\geq 0}\frac{1}{(m!)^2}z^{2m}.
  \end{equation}
\end{notation}
The following lemma explains why this function is relevant to us. It allows us to give an explicit expression for the partition function for the open necklace graph~$\graph_{1,2}$ from \cref{figure:necklace-graph}. The necklace graph is the dual graph of the two-holed torus as shown in \cref{figure:dualoftwoholedandnecklace}
\begin{figure}[t!]
  \centering
  \begin{subfigure}[b]{\textwidth}
    \centering
    \begin{tikzpicture}[scale=1.75]
      \node[vertex] (A) at (0,0) {};
      \node[vertex,fill] (B) at (1,0) {};
      \node (C) at (-1,0) {};
      \node (D) at (2,0) {};

      \draw (A) +(0,-0.1) node [below] {$1$};
      \draw (B) +(0,-0.1) node [below] {$2$};

      \draw (A) edge [bend left]  (B);
      \draw (A) edge [bend right] (B);
      \draw[middlearrow={>}] (A) to node [above] {$x$} (C);
      \draw[middlearrow={<}] (B) to node [above] {$y^{-1}$} (D);
    \end{tikzpicture}
    \caption{Open necklace graph}
    \label{figure:necklace-graph}
  \end{subfigure}

  \begin{subfigure}[b]{\textwidth}
    \centering
    \begin{tikzpicture}[scale=1.75]
      \node[vertex,fill] (A) at (0,0) {};
      \node[vertex,fill] (B) at (1,0) {};
      \node[vertex,fill] (C) at (2,0) {};
      \node[vertex,fill] (D) at (3,0) {};
      \node (E) at (3.5,0) {};
      \node (F) at (4,0) {};
      \node[vertex,fill] (G) at (4.5,0) {};
      \node[vertex,fill] (H) at (5.5,0) {};
      \node(I) at (-0.5,0){};
       \node(J) at (6.0,0){};

      \draw (A) edge [bend left]  (B);
      \draw (A) edge [bend right] (B);
      \draw (B) edge (C);
      \draw (C) edge [bend left]  (D);
      \draw (C) edge [bend right] (D);
      \draw (D) edge (E);
      \draw (F) edge (G);
      \draw (G) edge [bend left]  (H);
      \draw (G) edge [bend right] (H);
      \draw(I)edge (A);
       \draw(J)edge (H);

      \draw (3.75,0) node {$\ldots$};

\end{tikzpicture}
    \caption{Open necklace graph with many beads}
    \label{subfigure:closed-necklace-graph}
  \end{subfigure}

  \caption{Necklace graphs}
\end{figure}
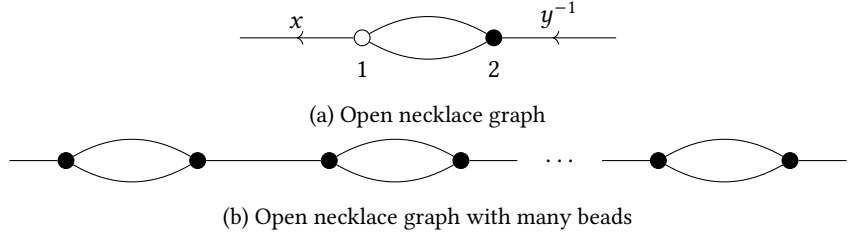

\begin{figure}[h]
	\centering
	\begin{tikzpicture}[baseline=(current bounding box.center)]
	\pic[
	scale=.5,
	rotate=90,
	every node/.style={transform shape},
	tqft,
	incoming boundary components=2,
	outgoing boundary components=1,
	offset=0.5,
	cobordism edge/.style={draw},
	every outgoing lower boundary component/.style={draw},
	every outgoing upper boundary component/.style={draw},
	name=a,
	];

	\pic[
	scale=.5,
	rotate=90,
	every node/.style={transform shape},
	tqft,
	incoming boundary components=1,
	outgoing boundary components=2,
	offset=-0.5,
	cobordism edge/.style={draw},
	anchor=outgoing boundary 1,
  every incoming lower boundary component/.style={draw, dashed},
	every incoming upper boundary component/.style={draw},
	name=b,
	];

	\draw (a-outgoing boundary 1) node [xshift=1em, yshift=-1em] {$\epsilon=1$};
	\draw (b-incoming boundary 1) node [xshift=-1em, yshift=-1em] {$\epsilon=0$};

	\node at (2.5, 0.5) {$=$};

	\node[vertex] (A) at (4,.5) {};
	\node[vertex,fill] (B) at (5,.5) {};
	\node (C) at (3,.5) {};
	\node (D) at (6,.5) {};

	\draw (A) +(0,-0.1) node [below] {$1$};
	\draw (B) +(0,-0.1) node [below] {$2$};

	\draw (A) edge [bend left]  (B);
	\draw (A) edge [bend right] (B);
	\draw[middlearrow={>}] (A) -- node [above] {$x$} (C);
	\draw[middlearrow={<}] (B) -- node [above] {$y^{-1}$} (D);
	\end{tikzpicture}
	\caption{Dual graph of the two-holed torus}
	\label{figure:dualoftwoholedandnecklace}
\end{figure}
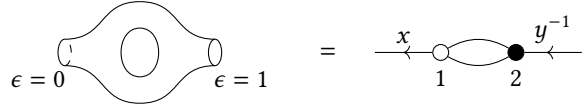
\begin{lemma}
  \label{lemma:necklace-bessel-product}
  Let~$\graph=\graph_{1,2}$ be the open necklace graph as in \cref{figure:necklace-graph}, with one half-edge oriented outwards and the other one oriented inwards. Then
  \begin{equation}
    \label{equation:necklace-bessel-product}
    \cK_{\graph,1}(t)=\bessel(t(x+y^{-1}))\bessel(t(x^{-1}+y)).
  \end{equation}
\end{lemma}

\begin{proof}
  The colored graph $(\graph,1)$ in the statement of the lemma is the dual graph obtained from a pair of pants decomposition of a two-holed torus with one hole oriented anticlockwise and the other one oriented clockwise. This is obtained by gluing two pairs of pants as shown in \cref{figure:cuttingtwoholes}.
	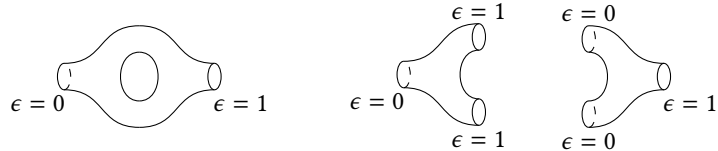
\begin{figure}[h]
		\centering
		\begin{tikzpicture}[baseline=(current bounding box.center)]
		\pic[
		scale=.5,
		rotate=90,
		every node/.style={transform shape},
		tqft,
		incoming boundary components=2,
		outgoing boundary components=1,
		offset=0.5,
		cobordism edge/.style={draw},
    every outgoing lower boundary component/.style={draw},
		every outgoing upper boundary component/.style={draw},
		name=a,
		];

		\pic[
		scale=.5,
		rotate=90,
		every node/.style={transform shape},
		tqft,
		incoming boundary components=1,
		outgoing boundary components=2,
		offset=-0.5,
		cobordism edge/.style={draw},
		anchor=outgoing boundary 1,
		every incoming lower boundary component/.style={draw, dashed},
		every incoming upper boundary component/.style={draw},
		name=b,
		];

		\draw (a-outgoing boundary 1) node [xshift=1em, yshift=-1em] {$\epsilon=1$};
		\draw (b-incoming boundary 1) node [xshift=-1em, yshift=-1em] {$\epsilon=0$};

		\pic[
		xshift=17em,
		scale=.5,
		rotate=90,
		every node/.style={transform shape},
		tqft,
		incoming boundary components=2,
		outgoing boundary components=1,
		offset=0.5,
		cobordism edge/.style={draw},
    every outgoing boundary component/.style={draw},
    every incoming lower boundary component/.style={draw, dashed},
	  every incoming upper boundary component/.style={draw},
		name=a,
		];

		\pic[
		xshift=10em,
		yshift=1.5em,
		scale=.5,
		rotate=90,
		every node/.style={transform shape},
		tqft,
		incoming boundary components=1,
		outgoing boundary components=2,
		offset=-0.5,
		cobordism edge/.style={draw},
    every outgoing boundary component/.style={draw},
    every incoming lower boundary component/.style={draw, dashed},
	  every incoming upper boundary component/.style={draw},
		name=b,
		];

		\draw (a-outgoing boundary 1) node [xshift=1em, yshift=-1em] {$\epsilon=1$};
		\draw (a-incoming boundary 2) node [yshift=1em] {$\epsilon=0$};
		\draw (a-incoming boundary 1) node [yshift=-1em] {$\epsilon=0$};
		\draw (b-incoming boundary 1) node [xshift=-1em, yshift=-1em] {$\epsilon=0$};
		\draw (b-outgoing boundary 2) node [yshift=1em] {$\epsilon=1$};
		\draw (b-outgoing boundary 1) node [yshift=-1em] {$\epsilon=1$};
		\end{tikzpicture}
		\caption{Pants decomposition for the two-holed torus}
		\label{figure:cuttingtwoholes}
	\end{figure}

	Now the graph $\graph_{1,2}$ is obtained two half-edges of a trivalent graph with one vertex as shown in \cref{figure:trivalentgluing}

	\begin{figure}[h]
		\centering
		\begin{tikzpicture}[baseline=(current bounding box.center)]
		\node[vertex] (A) at (0,.5) {};
		\node[vertex,fill] (B) at (1,.5) {};
		\node (C) at (-1,.5) {};
		\node (D) at (2,.5) {};

		\draw (A) edge [bend left]  (B);
		\draw (A) edge [bend right] (B);
		\draw[middlearrow={>}] (A) -- (C);
		\draw[middlearrow={<}] (B) -- (D);

		\node[vertex] (A) at (4,.5) {};
		\node[vertex,fill] (B) at (6,.5) {};
		\node (C) at (3,.5) {};
		\node (D) at (7,.5) {};

		\draw[middlearrow={<}, bend left] (A) to ++ (.8,0.5);
		\draw[middlearrow={<}, bend right] (A) to ++ (.8,-0.5);

		\draw[middlearrow={>}] (A) -- (C);
		\draw[middlearrow={<}] (B) -- (D);

		\draw[middlearrow={>}, bend right] (B) to ++ (-.8,0.5);
		\draw[middlearrow={>}, bend left] (B) to ++ (-.8,-0.5);
		\end{tikzpicture}
		\caption{Genus one graph from gluing }
		\label{figure:trivalentgluing}
	\end{figure}
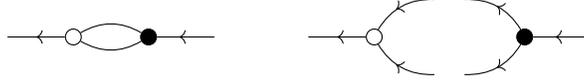
  By definition, \cref{proposition:cuttingandparitymatching}, and~\eqref{equation:formal-tqft}, we have
  \begin{align}
     \cK_{\graph, 1}(t)
     &=\left\langle \exp(t(\widetilde{W}_{v,0}(x,u^{-1},v^{-1})+\widetilde{W}_{v,0}(y^{-1},u,v)))\right\rangle_{u,v}\\
     &=\frac{1}{(2\pi{\sqrt{-1}})^2}\iint_{S^1\times S^1}\exp(t(\widetilde{W}_{v,0}(x,u^{-1},v^{-1})+\widetilde{W}_{v,0}(y^{-1},u^{-1},v^{-1})))\frac{\dd u}{u}\frac{\dd v}{v}\\
     &=\frac{1}{(2\pi{\sqrt{-1}})^2}\iint_{S^1\times S^1}\exp(t(\widetilde{W}_{v,0}(x,u,v)+\widetilde{W}_{v,0}(y^{-1},u,v)))\frac{\dd u}{u}\frac{\dd v}{v}\\
     \intertext{which we can interpret as}
     &=[\exp(t(\widetilde{W}_{v,0}(x,u,v)+\widetilde{W}_{v,0}(y^{-1},u,v)))]_{u^0v^0}.\\
     \intertext{Via the following sequence of rewrites}
     &=\sum_{\substack{a,b,c,d,a',b',c',d'\\b+d+b'+d'=a'+c'+a+c\\a'+d'+a+d=b+c+b'+c'}}\frac{{t}^{(a+b+c+d)+(a'+b'+c'+d')}}{a!a'!b!b'!c!c'!d!d'!}{x^{(a+b)-(c+d)}}y^{-(a'+b')+(c'+d')}\\
     &=\sum_{\substack{a,b,c,d,a',b',c',d'\\c+c'=d+d'\\a+a'=b+b'}}\frac{{t}^{(a+b+c+d)+(a'+b'+c'+d')}}{a!a'!b!b'!c!c'!d!d'!}x^{(a+b)-(c+d)}y^{-(a'+b')+(c'+d')}\\
     &=\sum_{\substack{a_1,b_1,c_1,d_1,a,a',c,c'\\ c_1+d_1=2(c+c')\\a_1+b_1=2(a+a')}}\frac{t^{a_1+b_1+c_1+d_1}}{a!(a_1-a)!a'!(b_1-a')!c!(c_1-c)!c'!(d_1-c')!}x^{a_1-c_1}y^{-b_1+d_1}\\
     &=\sum_{m,n}\sum_{\substack{a_1,b_1,c_1,d_1,a,c\\a\leq m, c\leq n\\c_1+d_1=2n;a_1+b_1=2m}}\frac{t^{2(m+n)}}{a!(a_1-a)!a'!(b_1-a')!c!(c_1-c)!c'!(d_1-c')!}x^{a_1-c_1}y^{-b_1+d_1}\\
     &=\sum_{m\geq 0,n\geq 0 }\sum_{\substack{a_1,b_1,c_1,d_1,a,c\\a\leq m, c\leq n\\c_1+d_1=2n; a_1+b_1=2m}}\frac{t^{2(m+n)}}{a_1!b_1!c_1!d_1!}x^{a_1-c_1}y^{-b_1+d_1} \binom{a_1}{a}\binom{b_1}{m-a}\binom{c_1}{c}\binom{d_1}{n-c}\\
     &=\sum_{m\geq 0,n\geq 0}\sum_{\substack{a_1,b_1,c_1,d_1\\ a_1+b_1=2m; c_1+d_1=2n}}t^{2(m+n)}\binom{2m}{m}\binom{2n}{n}\frac{x^{a_1-c_1}\left(y^{-1}\right)^{b_1-d_1}}{a_1!b_1!c_1!d_1!}\\
     \intertext{we finally obtain}
     &=\bessel(t(x+y^{-1}))\bessel(t(x^{-1}+y))
  \end{align}
  as desired.
\end{proof}

As a direct corollary of \cref{lemma:laplace-transform-interpretation} and \cref{remark:klugscheisser} we obtain
\begin{corollary}Let $\graph=\graph_{1,2}$ be the graph $\graph_{1,2}$ with no colored vertices, then
  \begin{equation}
    \cK_{\graph,0}(t)=\bessel(t(x+y))\bessel(t(x^{-1}+y^{-1})).
  \end{equation}
\end{corollary}

Before we discuss the general case, we will give a formula in the genus two case.
This formula will be revisited in \cite[Appendix~B]{gp-sympl}
\begin{corollary}
  \label{corollary:genus-2-epsilon-1}
  Let~$\graph$ be a genus two graph, without half-edges. We consider the case~$\epsilon=1$. Then the inverse Fourier--Laplace transform of the period~$\pi_{\widetilde{W}_{\graph,\epsilon}}(t)$ is given by
  \begin{equation}
    \sum_{n\geq 0}\frac{(2n!)^2}{n!^6}t^{2n}.
  \end{equation}
\end{corollary}

\begin{proof}
  If we cut the colored graph $(\graph,1)$ along any edge then we get back the graph considered in \cref{lemma:laplace-transform-interpretation}. Hence,
by \cref{lemma:laplace-transform-interpretation},~\eqref{equation:formal-tqft}, and \cref{lemma:necklace-bessel-product}, we have putting after~$x=y$
  \begin{equation}
    \begin{aligned}
      &\widecheck{\pi}_{\widetilde{W}_{\graph,1}}(t) \\
      &\quad=\hspace{-.4cm}\sum_{m\geq 0,n\geq 0}\hspace{-.4cm}\sum_{\substack{a_1,b_1,c_1,d_1\\ a_1+b_1=2m; c_1+d_1=2n}}\hspace{-.6cm}t^{2(m+n)}\binom{2m}{m}\binom{2n}{n}\frac{1}{a_1!b_1!c_1!d_1!}\left( \frac{1}{2\pi\sqrt{-1}}\int_{\sphere^1} x^{a_1-c_1}\cdot(x^{-1})^{b_1-d_1}\frac{\dd x}{x}\right) \\
      &\quad=\hspace{-.4cm}\sum_{m\geq 0,n\geq 0}\hspace{-.4cm}\sum_{\substack{a_1,b_1,c_1,d_1\\ a_1+b_1=2m; c_1+d_1=2n}}\hspace{-.6cm}t^{2(m+n)}\binom{2m}{m}\binom{2n}{n}\frac{1}{a_1!b_1!c_1!d_1!}\left( \frac{1}{2\pi\sqrt{-1}}\int_{\sphere^1}x^{a_1-c_1-b_1+d_1}\frac{\dd x}{x}\right)\\
      &\quad=[ \mathrm{B}(t(x+x^{-1}))\mathrm{B}(t(x^{-1}+x))]_{x^0}.
    \end{aligned}
  \end{equation}
  Using the definition of the twisted Bessel function, the Vandermonde identity and some elementary manipulations, we obtain
  \begin{equation}
    \begin{aligned}
      [\mathrm{B}(t(x+x^{-1}))^2]_{x^0}
      &=\sum_{n\geq 0}\left( \sum_{a+b=n}\frac{1}{a!^2b!^2} \right)[(x+x^{-1})^{2n}]_{x^0}t^{2n} \\
      &=\sum_{n\geq 0}\left( \sum_{a=0}^n\frac{n!^2}{a!^2(n-a)!^2} \right)\binom{2n}{n}\frac{t^{2n}}{n!^2} \\
      &=\sum_{n\geq 0}\binom{2n}{n}^2\frac{t^{2n}}{n!^2} \\
      &=\sum_{n\geq 0}\frac{(2n!)^2}{n!^6}t^{2n}.
    \end{aligned}
  \end{equation}
\end{proof}

\paragraph{Applying the machinery}
For~$g\geq 3$ we can describe an inductive procedure. Denote the product of Bessel functions~$\bessel(t(x+y)\bessel(t(x^{-1}+y^{-1})))$ by~$\rT_1(x,y)$.\index{T@$\rT_1(x,y)$}
Observe that
\begin{equation}
  \begin{aligned}
    \rT_1\left(x,y^{-1}\right)&=\rT_1\left(x^{-1},y\right)=\bessel(t(x+y^{-1})\bessel(t(x^{-1}+y)))\\
    \rT_1\left(x^{-1},y^{-1}\right)&=\rT_1\left(x,y\right)=\bessel(t(x+y)\bessel(t(x^{-1}+y^{-1}))).
  \end{aligned}
\end{equation}
This calculates the effect of changing the orientation of the boundary of the two-holed torus.
\begin{definition}
  Define inductively using convolution the function\index{T@$\rT_{k+1}(x,y)$}
  \begin{equation}
    \begin{aligned}
      \rT_{k+1}\left(x,y\right)&:=[\rT_k\left(x,z\right)\rT_1\left(z,y\right)]_{z^0}\\
\end{aligned}
  \end{equation}
\end{definition}

The following proposition is an application of the usual machinery of determining the partition function by cutting a closed surface to easier pieces, and it explains the definition of~$\rT_{k+1}\left(x,y\right)$.
\begin{proposition}
  \label{proposition:necklacegraph}
  Let~$\Sigma_{g,2}$ be a genus $g$ surface with two holes, one oriented anticlockwise and the other oriented clockwise. Then the partition function of~$\Sigma_{g,2}$ is given by $\rZ^{\operatorname{gp}}(t)(\Sigma_{g,2})=\rT_{g}\left(x,y^{-1}\right)$.
\end{proposition}

\begin{proof}
  Consider the necklace graph~$\graph_{1,2}$ as in \cref{lemma:necklace-bessel-product} show in \cref{figure:necklace-graph}. The genus two surface with two holes $\Sigma_{2,2}$ with one hole oriented anti-clockwise and the other hole oriented clockwise can be obtained gluing two $\Sigma_{1,2}$ as shown in \cref{figure:gluingtwotori}.

  In terms of the dual graph this is obtained by gluing two open necklace graphs \cref{figure:dualoftwoholedandnecklace} as shown in \cref{figure:dualgluingoftwobeadednecklace}.

\begin{figure}[h]
	\centering
  \begin{subfigure}[b]{\textwidth}
	  \begin{tikzpicture}[baseline=(current bounding box.center)]
	    \pic[
	    scale=.4,
	    rotate=90,
	    every node/.style={transform shape},
	    tqft,
	    incoming boundary components=2,
	    outgoing boundary components=1,
	    offset=0.5,
	    cobordism edge/.style={draw},
	    every outgoing lower boundary component/.style={draw},
	    every outgoing upper boundary component/.style={draw},
	    name=a,
	    ];

	    \pic[
	    scale=.4,
	    rotate=90,
	    every node/.style={transform shape},
	    tqft,
	    incoming boundary components=1,
	    outgoing boundary components=2,
	    offset=-0.5,
	    cobordism edge/.style={draw},
	    anchor=outgoing boundary 1,
      every incoming lower boundary component/.style={draw, dashed},
	    every incoming upper boundary component/.style={draw},
	    name=b,
	    ];

      \node at (1.75, 0.5) {$+$};

	    \pic[
	    xshift=10em,
	    scale=.4,
	    rotate=90,
	    every node/.style={transform shape},
	    tqft,
	    incoming boundary components=2,
	    outgoing boundary components=1,
	    offset=0.5,
	    cobordism edge/.style={draw},
	    every outgoing lower boundary component/.style={draw},
	    every outgoing upper boundary component/.style={draw},
	    name=c,
	    ];

	    \pic[
	    xshift=10em,
	    scale=.4,
	    rotate=90,
	    every node/.style={transform shape},
	    tqft,
	    incoming boundary components=1,
	    outgoing boundary components=2,
	    offset=-0.5,
	    cobordism edge/.style={draw},
	    anchor=outgoing boundary 1,
      every incoming lower boundary component/.style={draw, dashed},
	    every incoming upper boundary component/.style={draw},
	    name=d,
	    ];

	    \pic[
	    xshift=18em,
	    yshift=1.5em,
	    scale=.4,
	    rotate=90,
	    every node/.style={transform shape},
	    tqft,
	    incoming boundary components=1,
	    outgoing boundary components=2,
	    offset=-0.5,
	    cobordism edge/.style={draw},
      every incoming lower boundary component/.style={draw, dashed},
	    every incoming upper boundary component/.style={draw},
	    name=e,
	    ];

	    \pic[
	    scale=.4,
	    rotate=90,
	    every node/.style={transform shape},
	    tqft,
	    incoming boundary components=2,
	    outgoing boundary components=1,
	    offset=0.5,
	    cobordism edge/.style={draw},
	    at=(e-outgoing boundary 1),
	    name=f,
	    ];

	    \pic[
	    scale=.4,
	    rotate=90,
	    every node/.style={transform shape},
	    tqft,
	    at=(f-outgoing boundary 1),
	    incoming boundary components=1,
	    outgoing boundary components=2,
	    offset=-0.5,
	    cobordism edge/.style={draw},
	    name=g,
	    ];

	    \pic[
	    scale=.4,
	    rotate=90,
	    every node/.style={transform shape},
	    tqft,
	    incoming boundary components=2,
	    outgoing boundary components=1,
	    offset=0.5,
	    cobordism edge/.style={draw},
	    at=(g-outgoing boundary 1),
	    every outgoing lower boundary component/.style={draw},
	    every outgoing upper boundary component/.style={draw},
	    name=h,
	    ];

	    \node at (5.25, 0.5) {$=$};

	    \draw (a-outgoing boundary 1) node [xshift=1em, yshift=-1em] {$\epsilon=1$};
	    \draw (b-incoming boundary 1) node [xshift=-1em, yshift=-1em] {$\epsilon=0$};
	    \draw (c-outgoing boundary 1) node [xshift=1em, yshift=-1em] {$\epsilon=1$};
	    \draw (d-incoming boundary 1) node [xshift=-1em, yshift=-1em] {$\epsilon=0$};
	    \draw (e-incoming boundary 1) node [xshift=-1em, yshift=-1em] {$\epsilon=0$};
	    \draw (h-outgoing boundary 1) node [xshift=1em, yshift=-1em] {$\epsilon=1$};
	  \end{tikzpicture}
	  \caption{Gluing surfaces}
	  \label{figure:gluingtwotori}
  \end{subfigure}

  \begin{subfigure}[b]{\textwidth}
	  \begin{tikzpicture}
	    \node[vertex] (A) at (0,.5) {};
	    \node[vertex,fill] (B) at (1,.5) {};
	    \node (C) at (-1,.5) {};
	    \node (D) at (2,.5) {};

	    \draw (A) edge [bend left]  (B);
	    \draw (A) edge [bend right] (B);
	    \draw[middlearrow={>}] (A) -- node [above] {$x$} (C);
	    \draw[middlearrow={<}] (B) -- node [above] {$y^{-1}$} (D);

	    \node at (2.25, 0.5) {$+$};

	    \node[vertex] (A) at (3.5,.5) {};
	    \node[vertex,fill] (B) at (4.5,.5) {};
	    \node (C) at (2.5,.5) {};
	    \node (D) at (5.5,.5) {};

	    \draw (A) edge [bend left]  (B);
	    \draw (A) edge [bend right] (B);
	    \draw[middlearrow={>}] (A) -- node [above] {$y$} (C);
	    \draw[middlearrow={<}] (B) -- node [above] {$z^{-1}$} (D);

	    \node at (5.75, 0.5) {$=$};

	    \node[vertex] (A) at (7,.5) {};
	    \node[vertex,fill] (B) at (8,.5) {};
	    \node (C) at (6,.5) {};
	    \node[vertex] (D) at (9,.5) {};
	    \node[vertex,fill] (E) at (10,.5) {};
	    \node (F) at (11,.5) {};

	    \draw (A) edge [bend left]  (B);
	    \draw (A) edge [bend right] (B);
	    \draw (B) edge              (D);
	    \draw (D) edge [bend left]  (E);
	    \draw (D) edge [bend right] (E);
	    \draw[middlearrow={>}] (A) -- node [above] {$x$} (C);
	    \draw[middlearrow={<}] (E) -- node [above] {$z^{-1}$} (F);
	  \end{tikzpicture}
	  \caption{Dual picture}
	  \label{figure:dualgluingoftwobeadednecklace}
  \end{subfigure}
  \caption{Two-holed surfaces of genus two by gluing}
\end{figure}
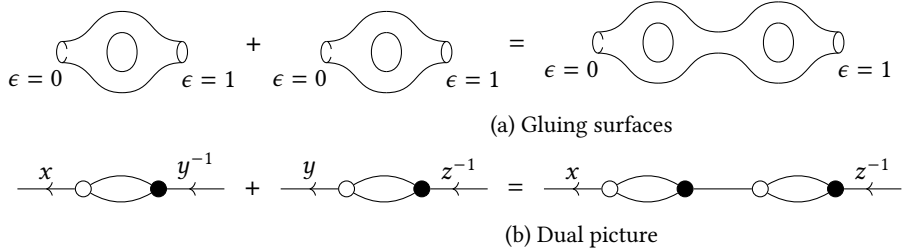

Then the partition function is given by
\begin{equation}
  \begin{aligned}
    \rZ^{\operatorname{gp}}(t)(\Sigma_{2,2})
    &=\langle \rT_1(x,z^{-1}),\rT_1(z,y^{-1})\rangle_{z}\\
    &=\frac{1}{2\pi\sqrt{-1}}\int_{\sphere^1}\rT_1(x,z^{-1})\rT_1(z^{-1},y^{-1})\frac{\dd z}{z}\\
    &=\rT_2(x,y^{-1}).
  \end{aligned}
\end{equation}
Hence by repeating this process, gluing copies of~$\graph_{1,2}$ to increase the genus, and applying the TQFT formalism given by \cref{proposition:tqft-from-wdvv} we obtained the required result.
\end{proof}

We now describe this result in terms of the period of graph potentials that we are interested in computing.
\begin{proposition}
  \label{proposition:steps}
  Let $(\graph_{g,2},\epsilon)$ be the (colored) dual graph of a surface with two holes obtained by gluing~$g$~copies of $\Sigma_{1,2}$ in a row. Then
  \begin{equation}
    \cK_{\graph_{g,2},\epsilon}(t)=\rT_{g}\left(x,y^{(-1)^{g}}\right)=\rT_{g}\left(x,y^{(-1)^{\epsilon}}\right).
  \end{equation}
\end{proposition}

\begin{proof}
  First of all observe that $\epsilon$ has the same parity as the genus $g$ of the surface by construction.

  Furthermore, if~$g$ is odd, then we need to make even number of cuts to get a disjoint union of $\Sigma_{1,2}$'s. Hence by \cref{proposition:cuttingandparitymatching}, we get that $\cK_{\graph_{g,2},\epsilon}(t)=\rZ^{\operatorname{gp}}(t)(\Sigma_{g,2})$. Then we are done by \cref{proposition:necklacegraph} when~$g$ is odd.

  Now if~$g$ is even, then consider the graph $\graph_{g-1,2}$ and an open necklace graph. Now by \cref{proposition:cutting-a-graph}, we get
  \begin{equation}
    \begin{aligned}
      \cK_{\graph_{g,2},0}(t)
      &=[\cK_{\graph_{g-1,2},1}(t),\rT_1\left(z,y^{-1}\right)]_{z^0}\\
      &=[\rT_{g-1}\left(x,z^{-1}\right),\rT_1\left(z^{-1},y\right)]_{z^0}\\
      &=\rT_g\left(x,y\right).
    \end{aligned}
  \end{equation}
\end{proof}
A direct very useful corollary of the above proposition that removes the restriction on the matching of the parity of coloring and genus in \cref{proposition:cuttingandparitymatching} is the following

\begin{corollary}
	\label{corollary:usefulsteps}
	Let $\graph_{g,2}$ be the open necklace graph of genus $g$ with two half-edges as shown in \cref{subfigure:closed-necklace-graph}. Let $c$ be a coloring of $\graph_{g,2}$ and $\epsilon$ be the parity of the coloring, then
  \begin{equation}
    \cK_{\graph_{g,2},c}(t)=\rT_{g}\left(x,y^{(-1)^{\epsilon}}\right).
  \end{equation}
\end{corollary}

\begin{proof}
  If the parity of~$c$ matches up with the parity of the genus as in \cref{proposition:cuttingandparitymatching}, then we are done. Otherwise choose one of the two half-edges and reverse the orientation of that half-edge. In the definition of $\cK_{\graph_{g,2},c}(t)$, we never integrate over a half-edge variable of the original graph $\graph$, hence we are now reduced to the situation in \cref{proposition:cuttingandparitymatching}.
\end{proof}

Now we can use \cref{proposition:steps} to get a formula for the periods of the \emph{closed} (i.e.~without half-edges) genus~$g$ trivalent colored graphs~$(\graph,c)$.

\begin{proposition}
  \label{proposition:partition-function}
  Let~$\graph$ be a trivalent graph of genus~$g\geq 2$ without leaves and~$c$ a coloring of the vertices. Let~$\epsilon\in\bF_2$ denote the parity of the coloring. Then
  \begin{equation}
    \cK_{\graph,c}(t)=\left[ \rT_{g-1}\left(x,x^{(-1)^{\epsilon}}\right) \right]_{x^0}.
  \end{equation}
\end{proposition}

\begin{proof}
  By \cref{corollary:potential-only-depends-on-parity} we know that~$\cK_{\graph,c}(t)$ only depends on the parity of the coloring. Hence we can assume that the graph has either one or zero colored vertices. Moreover since the periods are invariants under mutations, we can assume that $\graph$ is the closed necklace graph of genus $g$ as in \cref{figure:1closed-necklace-graph} with one or zero colored vertex.

  \begin{figure}[h]
  	\centering
  	\begin{tikzpicture}[scale=1.75]
  	\node[vertex] (A) at (0,0) {};
  	\node[vertex] (B) at (1,0) {};
  	\node[vertex] (C) at (2,0) {};
  	\node[vertex] (D) at (3,0) {};
  	\node (E) at (3.5,0) {};
  	\node (F) at (4,0) {};
  	\node[vertex] (G) at (4.5,0) {};
  	\node[vertex] (H) at (5.5,0) {};

  	\draw (A) edge [bend left]  (B);
  	\draw (A) edge [bend right] (B);
  	\draw (B) edge (C);
  	\draw (C) edge [bend left]  (D);
  	\draw (C) edge [bend right] (D);
  	\draw (D) edge (E);
  	\draw (F) edge (G);
  	\draw (G) edge [bend left]  (H);
  	\draw (G) edge [bend right] (H);

  	\draw (3.75,0) node {$\ldots$};

  	\draw[rounded corners] (H) -- (6,0) -- (6,-0.4) -- (-0.5,-0.4) -- (-0.5,0) -- (A);
  	\end{tikzpicture}
  	\caption{Closed necklace graph with many beads}
  	\label{figure:1closed-necklace-graph}
\end{figure}
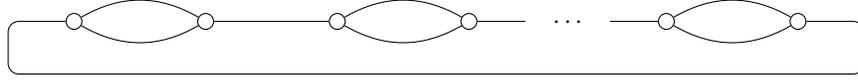

  First consider the case when $g$ and $c$ have the same parity. Then the result directly follows from \cref{proposition:steps} by applying \cref{proposition:cutting-a-graph}. Now in the case~$g$ and~$c$ have opposite parity, then to compute~$\cK_{(\graph_{g,0},\epsilon)}(t)$, cut an edge~$e$ of the graph $\graph$ to get the graph $(\graph_{g-1,2},\epsilon)$ considered in \cref{proposition:steps}. Then by we are done by first applying \cref{corollary:usefulsteps} and then \cref{proposition:cutting-a-graph}.

\end{proof}

\paragraph{The main result}
This proposition gives an effective method to compute the period of a graph potential, giving the main result of this section. It allows us to express (the inverse Fourier--Laplace transform of) the periods of the Laurent polynomials from \cref{section:graph-potentials} using trace-class operators in~$\LT$.

Define
\begin{itemize}
	\item $A$ to be the Hilbert--Schmidt operator given by~$\rT_1\left(x,y\right)$;
	\item $S$ to be the bounded linear operator~$S\left(\sum_{n\in\bZ}a_nx^n\right)=\sum_{n\in\bZ}a_nx^{-n}$.
\end{itemize}
We consider $A$ as matrix with respect to the orthonormal basis and reinterpret \cref{corollary:genus-2-epsilon-1} as follows.
\begin{lemma}
	\label{lemma:reinterpretinggenus-2-epsilon}
  Let $(\graph_2,1)$ be a genus two graph with one colored vertex, then
	\begin{equation}
	  \cK_{\graph_2,1}(t)=\tr(AS)=\sum_{n\geq 0}\frac{(2n!)^2}{n!^6}t^{2n}.
	\end{equation}
	In particular the composition $AS$ is trace class.
\end{lemma}

Since trace-class operators form an ideal, it follows that~$A^{a+1}S^{b}$ is also trace-class for non-negative integers~$a, b$. Moreover~$S$ commutes with~$A$.

The following is the important computational tool for computing periods of graph potentials.
\begin{theorem}
  \label{theorem:main-theorem-trace-class-operators}
  Let~$(\graph,c)$ be a colored trivalent graph of genus~$g\geq 2$ (without half-edges). Then,
  \begin{equation}
    \label{equation:laplace-transform-via-trace}
    \widecheck{\pi}_{\widetilde{W}_{\graph,c}}(t)=\tr(A^{g-1}S^{\epsilon+g})
  \end{equation}
  where~$\epsilon$ denotes the parity of the number of colored vertices in~$\graph$.
\end{theorem}

\begin{proof}
  By \cref{corollary:potential-only-depends-on-parity} we can assume that the number of colored vertices is either zero or one. To show the equality in~\eqref{equation:laplace-transform-via-trace} we use \cref{lemma:laplace-transform-interpretation} and \cref{proposition:partition-function}, so that we need to show that
  \begin{equation}
    [\rT_{g-1}(x,x^{(-1)^{\epsilon}})]_{x^0}=\tr(A^{g-1}S^{\epsilon+g}).
  \end{equation}
  We do this by induction. If~$g=2$, then this follows from \cref{corollary:genus-2-epsilon-1} for~$\epsilon=1$, whilst the case of~$\epsilon=0$ is analogous.

  Next, observe that by the symmetry properties of~$A$, we have that~$S$ commutes with~$A$.

  Now if~$g$ is arbitrary, we consider~$\rT_{g-1}(x,y)\in\LT^{\otimes2}$, writing it as~$\sum_{i,j\in\bZ}b_{i,j}x^iy^j$ for some matrix~$B=(b_{i,j})$ in the orthonormal basis~$\{x^i\}_{i\in \bZ}$. We claim that~$B=A^{g-1}S^{g-1}$. This follows by induction, using that~$S$ commutes with~$A$.

  Hence we get
  \begin{equation}
    \begin{aligned}
      [\rT_{g-1}(x,x^{(-1)^\epsilon})]_{x^0}
      &=\sum_{\substack{i,j\in\bZ \\ i+(-1)^{\epsilon} j=0}}b_{i,j} \\
      &=\sum_{\substack{i,j\in\bZ \\ i=(-1)^{\epsilon+1} j}}b_{i,j} \\
      &=\sum_{i\in\bZ}(S^{\epsilon+1}A^{g-1}S^{g-1})_{i,i} \\
      &=\tr(A^{g-1}S^{\epsilon+g}),
    \end{aligned}
  \end{equation}
  again using that~$S$ commutes with~$A$.
\end{proof}

\begin{remark}
  This gives an explicit and efficient method to compute the period sequences, by truncating the power series in~$t$. This method is independent of~$g$, and the number of periods that one can compute depends only on the degree of the truncation, and \emph{not} on~$g$.
  As an illustration of this procedure,
  we have collected some output in \cref{table:periods-odd,table:periods-even}.
  In \cite[Appendix~B]{gp-sympl} we will discuss some patterns in this table.

  Alternatively, using the main result of \cite{gp-sympl} one can compute quantum periods of these moduli spaces via the abelian/non-abelian correspondence in Gromov--Witten theory. But in this approach one fixes~$g$ and does a special analysis for each value of~$g$ (which is only feasible for low~$g$). Developing the details of the abelian/non-abelian correspondence in this case and comparing the two methods is left for future work.\end{remark}

\paragraph{Some connections to other works in the literature}
There exists a degeneration of our setup, which links it to mirror symmetry for Grassmannians of planes. For this we let~$\graph$ be a genus zero graph with~$n\geq 3$ half-edges and all vertices uncolored. The cardinality of the set of edges~$E$ is~$2n-3$ and we assign variables~$x_1,\ldots, x_{2n-3}$ as before to all the half-edges. Consider the substitution
\begin{equation}
  (x,y,z)\rightarrow \left(\frac{\tau}{X},\frac{Y}{\tau},\frac{Z}{\tau}\right).
\end{equation}
where $x,y,z$ are variables assigned to adjacent edges on a vertex~$v$; the $X,Y,Z$ are new variables and $\tau$ is a formal parameter. In this set-up
\begin{equation}
  \label{equation:NU}
  \begin{aligned}
    \widetilde{W}_{v,0}(x,y,z)&={xyz}+\frac{x}{yz}+\frac{y}{xz}+\frac{z}{xy}\\
    &=\tau^{-1} \left(\frac{YZ}{ X}+\tau^{4}\frac{1}{XYZ}+\frac{XY}{Z}+\frac{ZX}{Y}\right)
  \end{aligned}
\end{equation}
Now consider the Laurent polynomial~$\widetilde{W}_{\graph,0}=\sum_{v\in V} \widetilde{W}_{v,0}(x,y,z)$, and take~$\lim_{\tau \rightarrow 0}\tau \widetilde{W}_{\graph}$. This is a Laurent polynomial~$\cG_{\graph}$ in the variables~$X_1,\dots, X_{2n-3}$. This is exactly the Laurent polynomial considered by Nohara--Ueda in \cite[Theorem 1.6]{MR3211821}, and brings us to the following remark.

\begin{remark}
\label{remark:NUimp}
The above observation and the \tqft results
from \cref{theorem:tqft-introduction} and \cref{proposition:tqft-from-wdvv}
tell us that the potential functions from \textcite{MR3211821} associated to the Grassmannian of planes in an~$n$\dash dimensional vector space
also give a \tqft which one can call the {\em Grassmannian \tqft}.
To the best of our knowledge this was not observed before.

In particular, by the above discussion and \eqref{equation:NU}
the graph potential \tqft coming from moduli of bundles that we develop
recovers the Grassmannian \tqft as a limit.
This should be compared with a corresponding fact
that genus zero conformal blocks recover
invariants of tensor product representations in the limit.
\end{remark}

\begin{sidewaystable}
  \centering
  \begin{tabular}{cccccccccc}
    \toprule
    $g$ & $p_0$ & $p_2$ & $p_4$ & $p_6$ & $p_8$     & $p_{10}$   & $p_{12}$       & $p_{14}$        & $p_{16}$            \\ \midrule
    2   & 1     & 8     & 216   & 8000  & 343000    & 16003008   & 788889024      & 40424237568     & 2131746903000        \\ 3   & 1     & 0     & 384   & 23040 & 3265920   & 435456000  & 68263641600    & 11300889600000  & 1984905402480000     \\ 4   & 1     & 0     & 576   & 11520 & 8769600   & 1175731200 & 445839609600   & 115772770713600 & 41211916193448000    \\ 5   & 1     & 0     & 768   & 0     & 16853760  & 928972800  & 1378578432000  & 295708763750400 & 237075779068128000   \\ 6   & 1     & 0     & 960   & 0     & 27518400  & 232243200  & 3112327680000  & 299893321728000 & 795162277629720000   \\ 7   & 1     & 0     & 1152  & 0     & 40763520  & 0          & 5892216422400  & 133905855283200 & 2006716647119184000  \\ 8   & 1     & 0     & 1344  & 0     & 56589120  & 0          & 9963493478400  & 22317642547200  & 4248683870158728000  \\ 9   & 1     & 0     & 1536  & 0     & 74995200  & 0          & 15571407667200 & 0               & 7983708676751808000  \\ 10  & 1     & 0     & 1728  & 0     & 95981760  & 0          & 22961207808000 & 0               & 13760135544283128000 \\ \bottomrule
  \end{tabular}
  \caption{Period sequence for the odd graph potential}
  \label{table:periods-odd}

  \vspace{1cm}

  \centering
  \begin{tabular}{ccccccccccc}
    \toprule
    $g$ & $p_0$ & $p_2$ & $p_4$ & $p_6$ & $p_8$     & $p_{10}$   & $p_{12}$       & $p_{14}$        & $p_{16}$             & $p_{18}$ \\
    \midrule
    2   & 1     & 0     & 384   & 0     & 645120    & 0          & 1513881600     & 0               & 4132896768000        & 0 \\
    3   & 1     & 0     & 576   & 0     & 6350400   & 0          & 136604160000   & 0               & 3976941969000000     & 0 \\
    4   & 1     & 0     & 576   & 0     & 12640320  & 0          & 805929062400   & 0               & 80306439693480000    & 0 \\
    5   & 1     & 0     & 768   & 0     & 18144000  & 0          & 1915060224000  & 0               & 401643111149280000   & 0 \\
    6   & 1     & 0     & 960   & 0     & 27518400  & 0          & 3418888704000  & 0               & 1062973988196120000  & 0 \\
    7   & 1     & 0     & 1152  & 0     & 40763520  & 0          & 5953528627200  & 0               & 2211592605702480000  & 0 \\
    8   & 1     & 0     & 1344  & 0     & 56589120  & 0          & 9963493478400  & 0               & 4323671149117320000  & 0 \\
    9   & 1     & 0     & 1536  & 0     & 74995200  & 0          & 15571407667200 & 0               & 7994421145174464000  & 0 \\
    10  & 1     & 0     & 1728  & 0     & 95981760  & 0          & 22961207808000 & 0               & 13760135544283128000 & 0 \\
    \bottomrule
  \end{tabular}
  \caption{Period sequence for the even graph potential}
  \label{table:periods-even}
\end{sidewaystable}

 \section{Applications of graph potentials}
\label{section:applications}
As this is the first paper in a series
we will now outline the results in the next installments.

\subsection{Mirror symmetry and moduli of rank-2 bundles}
\label{subsection:mirror-symmetry}
Whilst they do not play a role in the current paper,
our original motivation to introduce graph potentials
was to study the algebro- and symplecto-geometric aspects of the moduli space~$\odd$
of rank~2 vector bundles with fixed determinant on a smooth projective curve~$C$.
The mirror dual of~$\odd$ is expected to be
a cluster-like variety equipped with a regular function,
the so-called Landau--Ginzburg potential (closely related to the Floer potential).
We propose in \cite{gp-sympl,gp-decomp}
that graph potentials can be seen as (building blocks of) mirrors to~$\odd$.

\paragraph{Enumerative mirror symmetry}
In \cite{gp-sympl} we will discuss aspects of the symplectic geometry of~$\odd$,
and relate it to the algebraic geometry of graph potentials.

In particular we will consider the \emph{quantum period}~$\mathrm{G}_{\odd}(t)$,
a generating function for closed genus zero Gromov--Witten invariants
with descendants and primary field being a point,
which can be also defined using the operator of quantum multiplication
by the first Chern class,
whose coefficients are closed genus zero correlators
without descendants but with two arbitrary cohomology insertions.

The invariant associated to graph potentials is the \emph{classical period},
the constant terms of powers of the Laurent polynomial.
The agreement of these two sequences of numbers is an important litmus test for mirror symmetry of Fano varieties,
and forms the main result of \cite{gp-sympl}.

Having established this,
the computational methods in \cref{section:tqft-computations} for~$\mathrm{Z}^{\mathrm{gp}}(t)$
will in turn allow to experiment with the quantum differential equation:
it provides the means to compute hundreds of coefficients.
This makes it possible to check properties of the quantum differential equation
in a highly non-trivial setting,
similar to what happens in \cite[\S9]{fano4folds}.

\paragraph{Homological mirror symmetry}
In \cite{gp-decomp} we will discuss certain decompositions that arise in the study of~$\odd$ and graph potentials.
Since $\odd$ is a Fano manifold,
its bounded derived category of coherent sheaves
is expected to have a homological mirror dual partner
in the sense of
Kontsevich's homological mirror symmetry conjecture, cf.
\cite{MR1403918} for Calabi--Yau varieties
and \cite[page~30]{Kontsevich-ENS} for the Fano version.

This mirror dual partner is expected to be a pair~$(Y,f)$,
where~$Y$ is a smooth quasiprojective variety and~$f$ a regular function on it.
The Laurent polynomials discussed before are then restrictions of the Landau--Ginzburg potential~$f$
to a torus~$\mathbb{G}_{\mathrm{m}}^n\subseteq Y$.
For more information on these pairs, and their ``tamings'' one is referred to \cite{MR3592695}.

Assuming that the categories of symplectic origin,
namely the so-called Fukaya category
of (graded) immersed (e.g.~embedded)
coisotropic (e.g.~Lagrangian)
subvarieties in $\odd$
decorated with a flat unitary connection
and the analogous Fukaya--Seidel category of
vanishing Lagrangian thimbles in $Y$, are well-defined,
the spaces $\odd$ and $(Y,f)$ are called \emph{homologically mirror dual}
if
\begin{itemize}
 \item the derived category of coherent sheaves on~$\odd$
            is equivalent to
      the derived Fukaya--Seidel category of~$(Y,f)$,
 \item the derived Fukaya category of~$\odd$
            is equivalent to
      the matrix factorization category of~$(Y,f)$.
\end{itemize}
We also refer the reader to
\cite[Conjecture~2.3]{MR2596637}
for a general discussion on Landau--Ginzburg models
and homological mirror symmetry conjectures.

A natural approach to tackle the first equivalence
starts with finding \emph{semiorthogonal decompositions} on either side of the mirror.
We propose in \cite[Conjecture~A]{gp-decomp} a conjectural semiorthogonal decomposition for $\dbcoh{\odd}$,
independently suggested by Narasimhan.
For a description of the state-of-the-art we refer to op.~cit.
As evidence for this conjecture we provide in \cite[Theorem~C]{gp-decomp} a \emph{motivic decomposition}
of~$\odd$ (resp.~$\dbcoh{\odd}$) in the Grothendieck rings of varieties and categories.

With respect to the graph potentials introduced in this paper,
it becomes interesting to study the second equivalence of categories.
Here one aims to find \emph{orthogonal decompositions} on either side of the mirror.
For the derived Fukaya category there is a natural decomposition in terms of eigenvalues of the quantum multiplication,
which is described by Mu\~noz \cite{MR1695800,MR1670396}.
In \cite[Theorem~B]{gp-decomp} we show that this eigenvalue decomposition is mirrored by a critical value decomposition of graph potentials.
This in turns gives further evidence for the conjectured semiorthogonal decomposition.

\paragraph{Reconstruction results}
In \cite{gp-torelli} we discuss how the Newton polytope of the graph potential
determines the graph and its coloring.
This is a reconstruction result
which on the algebro-geometric side mirror to graph potentials
corresponds to a \emph{combinatorial non-abelian Torelli theorem}.
We refer to op.~cit.~for more context and applications.

\subsection{Mirror approach to the Atiyah--Floer conjecture}
\label{subsection:atiyah-floer}
Finally, we speculate on a variation on the theme of the Atiyah--Floer conjecture,
suggested by mirror symmetry,
and how graph potentials could lead to progress.

The Atiyah--Floer conjecture from \cite{MR0974342} states that two homology theories,
both introduced by Floer,
are isomorphic.
One is the instanton Floer homology~$\operatorname{HF}^{\mathrm{inst}}(B)$ of a 3-manifold~$B$ \cite{MR0956166}.
The other is the symplectic Floer homology of two Lagrangian subvarieties~$R(B_-)$ and~$R(B_+)$
inside a symplectic variety~$R(\Sigma)$
associated to a Heegaard splitting~$B=B_-\cup_\Sigma B_+$ into two handlebodies along a common boundary surface~$\Sigma$ \cite{MR0965228}.
The notation~$R(*)$ stands for the moduli space of flat~$\mathrm{SU}(2)$-connections on~$*=\Sigma,B_\pm,B$.
The (smooth locus of the) variety~$R(\Sigma)$ is equipped with the Narasimhan--Atiyah--Bott--Goldman symplectic structure,
for which~$R(B_\pm)\to R(\Sigma)$ are Lagrangian embeddings.

Donaldson proposed to extend and categorify the Atiyah--Floer conjecture.
For an overview one is referred to \cite{MR3838878}.
Associated to the surface~$\Sigma$ one then has a category~$\mathcal{C}(\Sigma)$,
whilst the handlebodies~$B_\pm$ with boundary~$\partial B_\pm=\Sigma$ define objects in~$\mathcal{C}(\Sigma)$.
This assignment is subject to the condition that the morphism spaces are identified with the symplectic Floer homology.
This (provisional) extended 4-dimensional TQFT is known as \emph{Donaldson--Floer theory}.
The category~$\mathcal{C}(\Sigma)$ should be a Fukaya-like category for the symplectic variety~$R(\Sigma)$,
whereas the Lagrangian subvarieties~$R(B_\pm)$ define objects in it.

On the other side of the mirror we have graph potentials
(and more complicated Landau--Ginzburg models constructed using graph potentials),
and their categories of matrix factorizations.
Similar to how categories of matrix factorizations of the potential~$x^n$
appear as basic building blocks of Khovanov--Rozansky field theories,
the categories of matrix factorizations of graph potentials
are expected to be mirror dual to the Fukaya categories of the symplectic varieties~$R(\Sigma)$.

The upshot of this approach is that it is purely algebraic,
whereas on the symplectic side one needs to do complicated analysis.
Likewise, $R(\Sigma)$ is singular,
complicating the study of the Fukaya category even further.
In particular,
our work suggests how one could try to construct Donaldson--Floer theory on the other side of the mirror.

\renewcommand*{\bibfont}{\small}
\printbibliography

\emph{Pieter Belmans}, \url{pieter.belmans@uni.lu} \\
Department of Mathematics, Universit\'e de Luxembourg, 6, avenue de la Fonte, L-4364 Esch-sur-Alzette, Luxembourg

\emph{Sergey Galkin}, \url{sergey@puc-rio.br} \\
PUC-Rio, Departamento de Matem\'atica, Rua Marqu\^es de S\~ao Vicente 225, G\'avea, Rio de Janeiro, Brasil

\emph{Swarnava Mukhopadhyay}, \url{swarnava@math.tifr.res.in} \\
School of Mathematics, Tata Institute of Fundamental Research, 1 Homi Bhabha Road, Navy Nagar, Colaba, Mumbai 400005, India

\end{document}